\documentclass[11pt, reqno]{amsart}

\usepackage{amsmath}
\usepackage{amscd}
\usepackage{mathrsfs}
\usepackage{amsfonts}
\usepackage{amssymb}
\usepackage{enumerate}
\usepackage{setspace}
\usepackage{dsfont}

\usepackage{array}

\usepackage{color}

\usepackage[hyperindex, colorlinks=true]{hyperref}

\usepackage[margin=1in]{geometry}

\newcommand{\eqn}{\begin{eqnarray}}
\newcommand{\een}{\end{eqnarray}}

\newtheorem{theorem}{Theorem}[section]

\theoremstyle{definition}

\numberwithin{equation}{section}

\begin{document}
\title{Growth in the Muskat problem} 
\author[R. Granero-Belinch\'{o}n]{Rafael Granero-Belinch\'{o}n}
\address{Departamento  de  Matem\'aticas,  Estad\'istica  y  Computaci\'on,  Universidad  de Cantabria.  Avda.  Los  Castros  s/n,  Santander,  Spain.}
\email{rafael.granero@unican.es}
\author{Omar Lazar}
\address{Departamento de An\'alisis Matem\'atico \& IMUS, Universidad de Sevilla, C/ Tarifa s/n, Campus Reina Mercedes, 41012 Sevilla, Spain}
\email{omar.lazar@us.es}
\begin{abstract}
We review some recent results on the Muskat problem modelling multiphase flow in porous media. Furthermore, we prove a new regularity criteria in terms of some norms of the initial data in critical spaces ($\dot{W}^{1,\infty}$ and $\dot{H}^{3/2}$).
\end{abstract}
\maketitle

\section{Introduction}
The mathematical study of multiphase flow in porous media is a very active research area \cite{NB,bear,Musk}. Besides being mathematically challenging, this problem is also  physically interesting as it models oil extraction \cite{Muskat:porous-media, Muskat, buckley1941mechanism,hassanizadeh1990mechanics}, tumor growth \cite{F}, beach evolution \cite{thornton2014hele} or a geothermal reservoir \cite{CF}. The problem of studying the free boundary flow in porous media is also known as the Muskat problem \cite{Musk}. The purpose of this paper is to review some recent results on the Muskat problem and also to prove a new regularity criteria following the same approach as the one developed in \cite{CL}.

Flow (at relatively slow velocities) in porous media evolves according to Darcy's Law
\begin{subequations}\label{Darcy}
\begin{align}
\frac{\mu}{\kappa}u(x,y,t)+\nabla p(x,y,t)&=-\rho(x,y,t)G(0,1)^T,\text{ for }(x,y,t)\in\mathbb{R}^2\times[0,T]\\
\nabla\cdot u(x,y,t)&=0,\text{ for }(x,y,t)\in\mathbb{R}^2\times[0,T]\\
\partial_t\rho(x,y,t)+\nabla\cdot(u(x,y,t)\rho(x,y,t))&=0\text{ for }(x,y,t)\in\mathbb{R}^2\times[0,T],
\end{align}
\end{subequations}
where $p,u, \rho$ and $\mu$ are the pressure, velocity, density and viscosity of the incompressible fluids while $\kappa$ denotes the permeability of the medium and $G$ denotes the acceleration due to gravity. In what follows we will assume that, in an appropriate choice of units, $G=1$. Unless otherwise stated, we will also fix $\kappa=1$. Darcy's Law was derived heuristically by Henry Darcy in 1856 \cite{Darcy} (although it can be derived rigorously using homogenization techniques \cite{hornung1997homogenization,Tartar:incompressible-porous-medium-homogenization}). Remarkably, (\ref{Darcy}a) was derived independently by Hele-Shaw \cite{H-S,HeleShaw:motion-viscous-fluid-parallel-plates} when he was studying viscid flow between two parallel flat plates separated by a narrow distance. The mathematical literature on the (horizontal, \emph{i.e.} where gravity is neglected) Hele-Shaw cell problem is also large. The interested reader can refer to \cite{Peter, chen1993hele, elliott1982weak, escher1997classical, escher1998center, SCH, chang2016free} and the references therein.

Equation \eqref{Darcy} is a system of hyperbolic active scalar equations. Some other equation in this family are the famous surface quasi-geostrophic equation \cite{constantin1994singular,KiselevSQG,constantin1999formation,Omar,berselli2002vanishing,Chae2singularSQG,cor2}, the magnetogeostrophic equation \cite{Moffatt,Friedlander2,Friedlander3,Friedlander4,Friedlander5}, the Stokes system \cite{lemarie2016navier,bae2015global} or the 2D Euler equation in vorticity formulation \cite{majda2002vorticity,majda1996two}. The Muskat problem studies the particular type of solution where there are two different immiscible fluids, a fluid on top with label $+$ and a fluid below with label $-$, with properties given by $(\rho^+,\mu^+)$ and $(\rho^-,\mu^-)$ (or a fluid with $(\rho^-,\mu^-)$ and a dry zone with $\rho^+=\mu^+=0$) separated by a moving interface, parametrized as
\begin{equation}\label{Gamma}
\Gamma(t)=\{(x,y)\in\mathbb{R}^2, \;\;(x,y)=(z_1(\alpha,t),z_2(\alpha,t)),\;\alpha\in\mathbb{R}\},
\end{equation} 
for certain functions $z_i:\mathbb{R}\times \mathbb{R}^+\mapsto\mathbb{R}$. We observe that, in this paper, unless otherwise stated we will assume $\mu^+=\mu^-$.

Thus, the goal is that, by getting a smooth enough solution for the interface equation, we obtain a weak solution of the conservation law \eqref{Darcy} on the whole plane (of course, at the same time, one each phase, the restriction of the weak solution is a strong solution). Besides obtaining a solution, we would like to understand whether the solution exists for all time and the dynamical properties of this solution or, at the contrary, if the solution presents finite time singularities. Similar problems have been studied for other active scalars as the surface quasi-geostrophic equation by Rodrigo, Gancedo and Gancedo \& Strain \cite{RodrigoQG,gancedo2008existence,gancedo2014absence}.

Another motivation to study the problem \eqref{Darcy} comes from the fact that, in a certain sense, the Muskat problem is a sort of \emph{parabolic} version of the water waves problem. To observe that, we have to use Lagrangian coordinates. Indeed, in Lagrangian variables, system \eqref{Darcy} reads
\begin{align}\label{DarcyLag}
\frac{\mu}{\kappa}\frac{d}{dt}\eta+A^T\nabla q&=-\rho^-G(0,1)^T,
\end{align}
where $A=(\nabla\eta)^{-1}$, $q=p\circ\eta$ and $\eta$ is the Lagrangian coordinates. At the same time, the water waves problem can be written as
\begin{align}\label{WaterWavesLag}
\rho^-\frac{d^2}{dt^2}\eta+A^T\nabla q&=-\rho^-G(0,1)^T.
\end{align}
Thus, we see that the connection between the Muskat and the water waves problems resembles the link betwen the heat and the wave equations. Similarly, one can compare the asymptotic models for the Muskat problem in \cite{GScr} and for the water waves problem in \cite{CGSW}.

Using the previous parametrization for $\Gamma(t)$ \eqref{Gamma}, the two-phase Muskat problem with parameters $(\rho^+,1),$ $(\rho^-,1)$ is equivalent to the following nonlinear and nonlocal evolution system for the unknowns $z_i$
\begin{equation}\label{Muskat0}
\partial_t z(\alpha)=\frac{\bar{\rho}}{\pi}\int_\mathbb{R}\frac{z_1(\alpha)-z_1(\beta)}{|z(\alpha)-z(\beta)|^2}(\partial_\alpha z(\alpha)-\partial_\alpha z(\beta))d\beta,
\end{equation}
where the integral is understood in Cauchy principal value sense, \emph{i.e.}
$$
\int_\mathbb{R}=\lim_{\epsilon\rightarrow0}\int_{B(0,\epsilon)^c\cap B(0,\epsilon^{-1})},
$$
and $\bar{\rho}=\frac{\kappa(\rho^--\rho^+)}{2}$. We observe that every integral is taken in principal value sense from this point onwards. 

Equivalently, when the viscosities satisfy $\mu^+=\mu^-$ and the interface is assumed to be the graph of the function $f(x,t)$, the Muskat system \eqref{Darcy} (or analogously \eqref{Muskat0}) can be written as a single nonlocal, nonlinear equation for the interface \eqref{Gamma}:

\begin{equation}\label{Muskat1}
\aligned
&\partial_{t}f=  \frac{\bar{\rho}}{\pi}\ \partial_{x}\int_\mathbb{R} \arctan \left(\frac{f(x,t)-f(x-\alpha,t)}{\alpha}\right) \  d\alpha
\\
& f(0,x)=f_{0}(x).
\endaligned
\end{equation}

The previous formulation as a conservation law with a nonlocal and nonlinear flux was obtained in \cite{ccgs-10}. Remarkably, the Muskat problem can also be written in terms of oscillatory integrals as 
\begin{equation} \label{Muskat2}
\aligned
& \partial_ t f (t,x) = \frac{\bar{\rho}}{\pi} \int_\mathbb{R} \partial_{x}\left(\frac{f(x,t)-f(x-\alpha,t)}{\alpha}\right) \ \int_{0}^{\infty} e^{-\delta} \cos\left(\delta \left(\frac{f(x,t)-f(x-\alpha,t)}{\alpha}\right)\right) \ d\delta \ d\alpha \\
& f(0,x)=f_{0}(x).
\endaligned
\end{equation}
This latter formulation was observed in \cite{CL}.

\section{Notation and functional setting}
We denote
$$
\Delta_{\alpha} f \equiv\frac{f(x,t)-f(x-\alpha,t)}{\alpha}.
$$
Similarly,
$$
\delta_{y} f(x)=f(x)-f(x-y)\text{ and }\bar \delta_{y} f(x)=f(x)-f(x+y).
$$
We define the Calder\'on operator $\Lambda=\sqrt{-\Delta}$. On the Fourier side, this operator is given as the action of the multiplier $|\xi|$, \emph{i.e.}
$$
\widehat{\Lambda f}(\xi)=|\xi|\hat{f}(\xi).
$$
In an analogous manner, we consider the Hilbert transform $\mathcal{H}$ given as the action of the multiplier $-i\text{sgn}(\xi)$, \emph{i.e.}
$$
\widehat{\mathcal{H} f}(\xi)=-i\text{sgn}(\xi)\hat{f}(\xi).
$$
We shall use the homogeneous $L^2$-based Sobolev space  $\dot H^{s}$, $s \in\mathbb{R}^+$, which is endowed with the (semi)-norm
$$
\Vert f \Vert_{\dot H^{s}} = \Vert \Lambda^{s} f \Vert_{L^{2}}.
$$
We will also use the $L^p$-based Sobolev spaces, $W^{n,p}(\mathbb{R})$, which are defined as
$$
W^{n,p}=\left\{u\in L^p(\mathbb{R}), \partial_x^{n} u\in L^p(\mathbb{R})\right\},
$$
with (semi-)norm
$$
\|u\|_{\dot{W}^{n,p}}=\|\partial_x^n u\|_{L^p}.
$$

Similarly, we define the homogeneous Wiener spaces $\dot{A}^\alpha(\mathbb{R})$ as
\begin{equation}\label{Wienerinhomo}
\dot{A}^\alpha(\mathbb{R})=\left\{u(x)\in L^1(\mathbb{R}),\text{ such that } \|u\|_{\dot{A}^\alpha(\mathbb{R})}=\int_{\mathbb{R}} |\xi|^\alpha|\hat{u}(\xi)|d\xi\right\}.
\end{equation}

Let us recall the definition of the homogeneous Besov spaces  $\dot B^{s}_{p,q}(\mathbb R)$ (see \cite{B, RS, bahouri2011fourier}). Let $(p,q,s) \in [1,\infty]^2 \times \mathbb R$. Let $f$ be a tempered distribution  (which is such that its Fourier transform is integrable near 0), then the homogeneous Besov space $\dot B^{s}_{p,q}(\mathbb R)$ is the space endowed with the following (semi)-norm
$$
\Vert f \Vert_{\dot B^{s}_{p,q}} =   \left\Vert \frac{\Vert \mathds{1}_{ ]0,1[}(s) \delta_{y}f + \mathds{1}_{[1,2[}(s) (\delta_{y}+\bar \delta_{y} f)\Vert_{L^{p}}}{ \vert{y}\vert^{s} } \right\Vert_{L^{q}(\mathbb R, \vert y \vert^{-1} dy)}
$$

We shall use the following classical embeddings. Let $(p_1,p_2,r_1,r_2) \in [1,\infty]^4$, then
$$
\dot B^{s_1}_{p_1,r_1}(\mathbb R) \hookrightarrow \dot B^{s_2}_{p_2,r_2}(\mathbb R),
$$
where $s_1+\frac{1}{p_2}=s_2+\frac{1}{p_1}$ and $r_1\leq r_2$. We also have for all $(p_1, s_1) \in [2,\infty] \times \mathbb R$,
$$
\dot B^{s_1}_{p_1,r_1}(\mathbb R) \hookrightarrow \dot B^{s_1}_{p_1,r_2}(\mathbb R),
$$
for all $(r_1,r_2) \in \ ]1,\infty]$ such that $r_1 \leq r_2$. Let $(s_1,s_2) \in \mathbb R^2$ so that $s_1<s_2$, then for all $\theta \in ]0,1[$ and $(p,r) \in [1,\infty]^2$, we have the following interpolation inequality
\begin{equation} \label{interp}
\Vert f \Vert_{\dot B^{\theta s_1 +(1-\theta)s_2}_{p,1}} \leq \frac{C}{s_2-s_1} \left(\frac{1}{\theta}+\frac{1}{1-\theta}\right)\Vert f \Vert^{\theta}_{\dot B^{s_1}_{p,r}} \Vert f \Vert^{1-\theta}_{\dot B^{s_2}_{p,r}}.
\end{equation}
We shall use the following useful generalized Calder\'on commutator type estimate (see e.g.   Dawson, McGahagan, and Ponce \cite{DMP} for a proof). Let $\Phi \in \dot W^{k+l,\infty}$ and let us consider the commutator  
$$
\left[\mathcal{H},\Phi \right] f=\mathcal{H}(\Phi f)-\Phi\mathcal{H}f.
$$
Then, for all $p\in]1,\infty[$ and $(k,l) \in \mathbb{N}$ 
\begin{equation} \label{cz}
\left\Vert \left[\mathcal{H},\Phi \right]  \partial^{k}_{x} f \right\Vert_{\dot W^{l,p}} \leq C_{k,l} \Vert \Phi \Vert_{\dot W^{k+l,\infty}} \Vert f \Vert_{L^{p}},
\end{equation}
for all $f \in L^{p}$. \\

Throughout the article,  $A \leq B$  means that there exists a constant $C>0$ depending only on controlled quantities such that $A \leq C B$. 

\section{Well-posedness}
Linearizing \eqref{Muskat0}, we obtain that the linear problem is
$$
\partial_t f=-\bar{\rho} \partial_x \mathcal{H} f.
$$
We note that the sign of $\bar{\rho}$ is crucial in the evolution. When $\bar{\rho}>0$, the linear problem reduces to a (fractional) heat equation, and it is therefore trivially well-posed in Sobolev spaces. However, when $\bar{\rho}<0$, the linear problem has an anti-diffusive character that makes the problem ill-posed in Sobolev spaces but well posed for analytic functions. The condition on the sign of $\bar{\rho}$ states that the fluids are in the stable regime if the lighter fluid is above the heavier fluid. Hence, it is a condition on the stratification of the fluids.

Due to (\ref{Darcy}a), the condition on the sign of $\bar{\rho}$ is equivalent to
\begin{equation*}
RT(t)=-(\nabla p^-(\Gamma(t))-\nabla p^+(\Gamma(t)))\cdot n>0,
\end{equation*}
where $n$ denotes the (upward) normal to $\Gamma(t)$. This latter condition is the well-known Rayleigh-Taylor stability condition \cite{Rayleigh:instability-jets,SaffTay}. This stability condition is ubiquitous in free boundary problems and it appears also when studying the water waves problem or the free-surface Euler equation \cite{cordoba2009rayleigh, coutand2007well}.

This stability condition that appears when studying the linear problem has to be taken into account when dealing with the full nonlinear problem \eqref{Muskat1} (unless surface tension effects are considered). 

Recalling \eqref{Muskat0} for $(x,f(x,t))$, we observe that the equation is invariant by the scaling
\begin{equation}\label{scaling}
f_\lambda(x,t)=\lambda^{-1}f(\lambda x,\lambda t).
\end{equation}
We observe that there are several spaces whose norm is also left invariant by this scaling. These spaces are called \emph{critical} for this equation. Three examples of critical spaces are
$$
L^\infty(0,T;\dot{W}^{1,\infty}), \;\; L^\infty(0,T;\dot{H}^{3/2}), \text{ and } L^\infty(0,T;\dot{B}^{1}_{\infty,\infty}).
$$
Spaces with more regularity (as for instance $L^\infty(0,T;\dot{H}^{2})$) are called subcritical while spaces with less regularity (as for instance $L^\infty(0,T;\dot{H}^{1})$) are called supercritical. The heuristic idea is that it is \emph{easy} to construct solutions in subcritical spaces and \emph{very difficult} to construct solutions in supercritical spaces. 

Due to all this, we see that there are three main ingredients we have to take into account when proving a well-posedness result:
\begin{enumerate}
\item the fluids need to have the good stratification (\emph{i.e.} the heavy fluid has to lie below the lighter fluid),
\item we should be able to parametrize the interface $\Gamma(t)$ as the graph of certain function $f(x,t)$ (otherwise, there exists some region where the fluid have the bad stratification),
\item the function $f$ needs to have subcritical (or, at most, critical) regularity.
\end{enumerate}

\subsection{Local existence}
Following the previous discussion, the \emph{basic} local existence result reads as follows
\begin{theorem}\label{teo1}Let $\rho^+<\rho^-$ be two fixed parameters. Fix $s=3$. Assume that $f_0\in \dot{H}^s(\mathbb{R})\cap L^2$. Then, there exists $0<T=T(\|f_0\|_{L^2},\|f_0\|_{\dot{H}^s})$ and a unique solution to \eqref{Muskat1}
$$
f\in C([0,T],L^2\cap\dot{H}^s).
$$
Furthermore, if $T=T_{max}<\infty$, then
$$
\limsup_{t\rightarrow T_{max}}\|f(t)\|_{C^{2+\delta}}=\infty.
$$
\end{theorem}

Physically, the previous hypotheses mean that the internal wave has no turning points (it is given as a graph) and separates fluids having the good stratification. More geometrically, the required smoothness on the data precludes consideration of initial data with a cusp. Actually, as noted in \cite{ambrose2004well}, (at the time) \emph{it is(was) an open problem as to whether the problem is well-posed for initial data with a cusp.}

Theorem \ref{teo1} was proved by C\'ordoba \& Gancedo \cite{c-g07} using energy methods (a similar result for the case where the spatial domain is a strip, for the case of a porous medium with two different permeabilities and for the case of three fluids was proved by C\'ordoba, Granero-Belinch\'on and Orive \cite{CGO} and Berselli, C\'ordoba and Granero-Belinch\'on \cite{BCG} and C\'ordoba \& Gancedo \cite{cordoba2010absence}, respectively). A different proof (using a formulation for \eqref{Muskat1} based on the tangent angle and arclength) was given by Ambrose \cite{ambrose2004well, AmbroseST} (see also \cite{tofts2017existence}). Another approach is the one by Escher, Matioc \& Walker \cite{emw15} where the authors used semigroup theory to obtain the similar result when the initial data is in the \emph{little  H\"older} spaces $h^{2+\delta}\subset C^{2+\delta}$ (see also the papers by Escher \& Matioc \cite{e-m10} and  Escher, Matioc \& Matioc \cite{escher2011generalized} where the case with small initial data is studied).

In order the initial regularity $\dot{H}^s$ can be relaxed for $3/2<s\leq5/2$ (so it allows for interfaces whose curvature is not bounded pointwise) some new ideas were required. Mathematically, an $H^s$ well-posedness result is challenging because the \emph{standard} energy estimates suggest $\|h\|_{C^{2+\delta}}$ as the quantity one needs to control. In the case of a fluid and a dry zone, \emph{i.e.} where the upper fluid is replaced by a dry zone, Cheng, Granero-Belinch\'on \& Shkoller \cite{CGS} introduced a new method to analyze \eqref{Darcy} and prove the local existence of an $\dot{H}^2$ solution. As the domain in \eqref{Darcy} $\Omega(t)$ is unknown, these authors first pull-back \eqref{Darcy} onto a fixed-in-time reference domain. By doing this, \eqref{Darcy} is transformed into a system of equations set on a fixed reference domain $\Omega$, but having time-dependent coefficients. Then, this new method is based on the analysis of the resulting quasilinear system of partial differential equations and combines new energy estimates in the bulk of the fluid with estimates for the interface (see also Shkoller \& Granero-Belinch\'on \cite{GS} for the case of two permeabilities). Due to the fact that this approach does not rely on the explicit structure of the singular integral equation \eqref{Muskat0}, it can be applied to study general domain geometries and permeability functions. 

In the case of two fluids with same viscosity, Constantin, Gancedo, Shvydkoy \& Vicol \cite{CGSVfiniteslope} proved the local existence of solution for $W^{2,p}$ $p\in(1,\infty]$ initial data. The proofs exploit the nonlocal nonlinear parabolic nature of \eqref{Muskat1} through a series of nonlinear lower
bounds for the nonlocal operators involved. Furthermore, these authors also prove that, as long as the slope of the interface remains uniformly bounded, the curvature remains bounded. A related result is the one by Pr\"{u}ss and Simmonet \cite{pruess2016muskat} where the authors prove the local existence for small initial data in $W^{2+\frac{1}{p},p}$.

Matioc \cite{matioc2016muskat}, by rewitting the Muskat problem as an abstract evolution equation in an appropriate functional setting, was able to prove the local existence for arbitrary $H^s$, $3/2<s<2$ initial data (see also \cite{matioc2017well,MR3841857,matioc2018well}).

For the case of two different viscosities in the RT stable regime the results are more scarce: C\'ordoba, 
C\'ordoba \& Gancedo proved the local existence for $H^3$ curves \cite{c-c-g10} and $H^4$ surfaces \cite{Cordoba-Cordoba-Gancedo:muskat-3d} (see also \cite{pernas2017local} for the case of two different permeabilities) while Cheng, Granero-Belinch\'on \& Shkoller prove the results for a $H^2$ graph.

\subsection{Global existence}
Before we can go over global existence results, we need to identify quantities that can be bounded for all positive times. The first of such results appeared in \cite{c-g09} and establishes the decay of $\|f(t)\|_{L^\infty}$:
\begin{theorem}\label{teo2}Let $\rho^+<\rho^-$ be two fixed parameters. Then the solution to the Muskat problem \eqref{Muskat1} satisfies
$$
\|f(t)\|_{L^\infty}\leq \|f_0\|_{L^\infty}
$$
\end{theorem}
To prove this result, C\'ordoba \& Gancedo used a pointwise estimate to compute the evolution of
$$
M(t)=\max_{x} f(x,t).
$$ 
For a similar result for the case of a bounded porous medium we refer to \cite{CGO}. Physically, this theorem means that the amplitude of the internal wave in a porous medium decays. 

Similarly, one can also prove an $L^2$ energy balance
\begin{theorem}\label{teo3}Let $\rho^+<\rho^-$ be two fixed parameters. Then the solution to the Muskat problem \eqref{Muskat1} satisfies
$$
\|f(t)\|_{L^2(\mathbb{R})}^2+\int_0^t\|u(s)\|_{L^2(\mathbb{R}^2)}^2ds\leq \|f_0\|_{L^2(\mathbb{R})}^2,
$$
or, equivalently,
$$
\|f(t)\|_{L^2(\mathbb{R})}^2+\frac{\rho^--\rho^+}{2\pi}\int_0^t\int_\mathbb{R}\int_\mathbb{R}\log\left(1+\left(\frac{f(x,s)-f(y,s)}{x-y}\right)^2\right) dxdyds\leq \|f_0\|_{L^2(\mathbb{R})}^2.
$$
\end{theorem}
The previous result was given by Constantin, C\'ordoba, Gancedo \& Strain \cite{ccgs-10} (see also \cite{BCG, CGS}).

Although the solution enjoys this decay of the \emph{relatively strong} $L^\infty$ norm and a energy balance that controls the velocity in $L^2(0,\infty;L^2(\mathbb{R}^2))$, this is not enough to obtain a global existence result of any kind. Then we have to turn our attention to other \emph{stronger} norms.

In that regards, the first result was given by C\'ordoba \& Gancedo \cite{c-g09}, where these authors proved that
\begin{theorem}\label{teo2b}Let $\rho^+<\rho^-$ be two fixed parameters and assume that
$$
\|f_0\|_{\dot{W}^{1,\infty}}<1. 
$$
Then the solution to the Muskat problem \eqref{Muskat1} satisfies
$$
\|f(t)\|_{\dot{W}^{1,\infty}}\leq \|f_0\|_{\dot{W}^{1,\infty}}.
$$
\end{theorem}
This result gives conditions ensuring the decay of a \emph{critical} norm. In the same spirit one has the analog result in terms of Wiener spaces (see \cite{ccgs-10,ccgs-13} and also the related work \cite{patel2017large})
\begin{theorem}\label{teo2c}Let $\rho^+<\rho^-$ be two fixed parameters and assume that
$$
\|f_0\|_{\dot{A}^{1}}<1/3. 
$$
Then the solution to the Muskat problem \eqref{Muskat1} satisfies
$$
\|f(t)\|_{\dot{A}^{1}}\leq \|f_0\|_{\dot{A}^{1}}.
$$
\end{theorem}

Theorems \ref{teo2b} and \ref{teo2c} imply that the internal wave will not break if the initial slope satisfies certain size restrictions (when measured in appropriate critical norms). \\

We have to distinguish two different global existence results: \\

\begin{enumerate}
\item global \emph{weak} solutions (when the equation \eqref{Muskat1} is satisfied in distributional sense)
\item global \emph{classical} solutions (when the equation \eqref{Muskat1} holds pointwise). \\
\end{enumerate}

Using Theorem \eqref{teo2b}, Constantin, C\'ordoba, Gancedo \& Strain \cite{ccgs-10} proved the global existence of weak Lipschitz solutions. These result was later extended to the case where the porous medium is bounded by one of the authors \cite{G}. Remarkably, in order for the wave to not break-down in the case where the domain is bounded, not only the slope has to be suitably small, but also the amplitude, the depth and the slope of the internal wave have to satisfy appropriate (explicit) conditions. Roughly speaking, these extra conditions linking the depth, the amplitude and the slope of the wave mean that the amplitude can not be such that the wave is \emph{close} to the bottom and that, the bigger the amplitude is, the smaller the slope has to be. 

Using Theorem \eqref{teo2c}, Constantin, C\'ordoba, Gancedo \& Strain \cite{ccgs-10} and Constantin, C\'ordoba, Gancedo, Rodr\'iguez-Piazza \& Strain \cite{ccgs-13} proved the global existence of classical solutions for initial data satisfying $\|f_0\|_{\dot{A}^1}<0.2$. The early work by C\'ordoba \& Gancedo \cite{c-g07} has already a global existence result for \eqref{Muskat1} in the spirit of Theorem \ref{teo2c}, but with a non-explicit size restriction $\|f_0\|_{\dot{A}^1}<\epsilon$ for \emph{certain} $\epsilon$.

With an initial data with suitably small Lipschitz norm, Constantin, Gancedo, Shvydkoy \& Vicol \cite{CGSVfiniteslope} proved the global existence of classical solution. This result, by further restricting the size of $\|f_0\|_{\dot{W}^{1,\infty}}$ establishes that the weak solution (which exists due to \cite{ccgs-10}) is, indeed, a classical solution. This result was later extended to the full range $\|f_0\|_{\dot{W}^{1,\infty}}<1$ by Cameron \cite{cameron2017global} (actually, the result by Cameron  is more general as the criteria is given in terms of the product of the supremum and infimum of the slope).\\

Another global existence result in a critical space (which allows the slope to be arbitrarily large) is the one by C\'ordoba \& Lazar \cite{CL}:
\begin{theorem}\label{teo4}Let $\rho^+<\rho^-$ be two fixed parameters and assume that $f_0\in L^2\cap\dot{H}^{3/2}\cap \dot{H}^{5/2}$ such that $\|f_0\|_{\dot{H}^{3/2}}$ is suitably small. Then, there exists a unique global strong solution of \eqref{Muskat2}
$$
f\in L^\infty(0,T;H^{5/2})\cap L^2(0,T;H^3),\;\forall\,T>0. 
$$
\end{theorem}

Some other global existence results are those in \cite{CGS} (smallness of the initial data in $H^2$, see also \cite{matioc2017well,matioc2018well}), \cite{gancedo2017muskat} (smallness in $\dot{A}^1$ but allows for two different viscosities) and \cite{escher2011generalized} (smallness in the little H\"older space $h^{2+\delta}$, see also \cite{e-m10}).

\section{Finite time singularities}
The first type of singularity for the Muskat problem is the \emph{turning singularity} and was proved by Castro, C\'ordoba, Fefferman, Gancedo \& L\'opez-Fern\'andez \cite{ccfgl} (see also \cite{CGO, BCG}). In these singularities, the initial interface is assumed to be a smooth graph and then, in a finite time $T_{breaking}$,
$$
\limsup_{t\rightarrow T_{breaking}}\|f(t)\|_{\dot{W}^{1,\infty}}=\infty.
$$
This singularity implies that the internal wave breaks, \emph{i.e.} after time $T_{breaking}$, the internal wave cannot be parametrized as a graph. Equivalently, this singularity means that the Muskat problem leaves the RT stable regime in finite time. However, a turning singularity do not imply a loss in derivatives for the solution.

The precise statement is
\begin{theorem}\label{teo5}Let $\rho^+<\rho^-$ be two fixed parameters. Then, there exist smooth initial data such that the unique strong solution of \eqref{Muskat1} satisfies
$$
\limsup_{t\rightarrow T_{breaking}}\|f(t)\|_{\dot{W}^{1,\infty}}=\infty,
$$
for $0<T_{breaking}<\infty$.
\end{theorem}

The proof of this result has several steps:
\begin{enumerate}
\item First, one considers the Muskat problem in its formulation for arbitrary curves \eqref{Muskat0} and proves the local existence of solution (forward and backward in time) for initial data who are analytic via a Cauchy-Kovalevsky theorem.
\item Then, one identifies $\partial_{\alpha}v_1(z(\alpha,t),t)$ as the quantity to track. Indeed, for initial data who are \emph{'about to break'} such that
$$
\partial_\alpha z(\alpha,t)\bigg{|}_{\alpha=0}=(0,1),
$$
\emph{i.e.} whose tangent vector is vertical, one has that, the sign condition
$$
\partial_{\alpha}v_1(z(\alpha,0),0)\bigg{|}_{\alpha=0}<0
$$  
is equivalent to breaking.
\item Now one constructs initial data such that $\partial_{\alpha}v_1(z(0,0),0)<0$ and $\partial_{\alpha}z(0,0)=(0,1)$.
\item Finally, one takes an analytic initial data such that the previous condition holds (using some mollification argument) and invokes the Cauchy-Kovalevsky Theorem in step 1. Using the forward and backward existence, we conclude the existence of internal waves (that can be parametrized as a graph for $-\delta<t<0$ that reach the initial data constructed in step 3 at time $t=0$ and that cannot be parametrized as a graph for $0<t<\delta$.
\end{enumerate}
We would like to remark that the analytic curves remain valid classical solutions to the Muskat problem \eqref{Muskat0} for $0<t<\delta$.

It is also interesting to note that the initial data leading to breaking waves can be taken with arbitrary small amplitude \cite{granero2013inhomogeneous}. 

G\'omez-Serrano \& Granero-Belinch\'on \cite{GG} (see also \cite{cordoba2015note,cordoba2017note}), by using a computer assisted proof together with a variation of the previous ideas, were able to study the effect of finite depth and varying permeability. Among other results, these authors showed that the existence of top and bottom for the porous medium can enhance the formation of turning singularities in the sense that there exists initial curves such that, when the depth is finite, they break in finite time, but, if the depth is infinite, they become smooth graphs.

In terms of loss of derivatives, Castro, C\'ordoba, Fefferman \& Gancedo \cite{castro2012breakdown} proved that there exist analytic initial data in the RT stable regime for the Muskat problem such that the solution turns to the RT unstable regime and later breaks down \emph{i.e.} no longer belongs to $C^4$.

Finally, another possible singularity is the self-intersection of the interface or the intersection of two different interfaces (for the problem with three different fluids or two fluids and a dry zone) while the curve remain itself smooth. These singularities are known as \emph{splash singularities} (when the self intersection happens at a single point) or \emph{splat singularities} (when the self-intersection happens along an interval). In this case, it was proved by Castro, C\'ordoba, Fefferman \& Gancedo \cite{ccfgonephase} that splash singularities can indeed occur for the one-phase Muskat problem, while the case of splat singularities was disregarded by C\'ordoba \& Pernas-Casta\~no \cite{cponephase}. When the case of several fluids is considered, Gancedo \& Strain \cite{gancedo2014absence} proved that splash/splat singularities cannot occur in finite time (see for the analog result for the case of the Euler equations \cite{FIL,coutand2016impossibility}).

\section{Wild solutions and mixing}
The Muskat problem \eqref{Muskat1} is \emph{ill-posed} in the RT unstable regime \cite{c-g07, CGO}. Thus, the construction of weak solutions to \eqref{Darcy} via the construction of classical solutions to \eqref{Muskat1} fails in the RT unstable regime. This leaves open the existence of weak solutions in the RT unstable regime and its dynamical properties.

Another topic that has attracted a lot of interest recently in the mathematical community is the construction of \emph{wild solutions} to different fluid dynamical problems. These wild solutions are weak solutions that have compact support in space and time and, thus, they break the uniqueness. In the Muskat problem these solutions are particularly interesting since they can be related to mixing of the fluids in the RT unstable regime. 

Thus, the existence of (possibly infinitely many) weak solutions in the RT unstable regime and their link to the mixing of the fluids appears as a very interesting research topic.

In this regards, it was proved by Sz{\'e}kelyhidi Jr \cite{szekelyhidi2012relaxation} the existence of weak solutions in the RT unstable regime (see also the works by F{\"o}rster \& Sz{\'e}kelyhidi Jr.\cite{forster2017piecewise} and Otto \cite{otto1999evolution,otto2001evolution}). Also, Castro, C\'ordoba \& Faraco \cite{castro2016mixing} proved that, starting with a smooth interface in the RT unstable regime, there exists a weak solution such that a \emph{mixing strip-like} region opens around the interface (see also Castro, Faraco \& Mengual \cite{castro2018degraded}). In other words, the free boundary assumption is replaced by the opening of this mixing zone where the fluids begin to mix.

\section{A new result} \label{h12}
In this section, we shall prove the following theorem
\begin{theorem} \label{h1/2}
Let $T>0$, assume that $f_{0} \in \dot H^{1/2} \cap \dot W^{1,\infty}$ then, if  
\begin{itemize}
 \item the $K=L^{\infty}([0,T],\dot W^{1,\infty})$ norm of the corresponding solution remains bounded
 \item $\Vert f_{0} \Vert_{\dot H^{3/2}}< C(K)$ is preserved
 \end{itemize}
 Then, the solution is global in time and we have
\begin{equation*}
 \Vert f \Vert^{2}_{\dot H^{1/2}}(T) + \frac{\pi}{1+K^2} \int_{0}^T \Vert f \Vert^{2}_{\dot H^{1}} \ ds \lesssim \Vert f_{0} \Vert^2_{\dot H^{1/2}} +  P\left(\Vert f \Vert_{L^{\infty}([0,T],\dot H^{3/2})}\right) \int_{0}^T \Vert f \Vert^{2}_{\dot H^{1}} \ ds
\end{equation*}
where $P(X)=X+X^2$
\end{theorem}
\noindent {\bf{Proof of Theorem \ref{h1/2}}}
For the sake of notational simplicity, we take $\bar{\rho}=\pi$. We do $\dot H^{1/2}$ estimates. Using \eqref{Muskat2}, we have that
\begin{eqnarray*}
\frac{1}{2}  \partial_{t} \Vert f \Vert^{2}_{\dot H^{1/2}}  &=& \int \Lambda^{1/2}f   \ \int_{0}^{\infty}   e^{-\delta}  \Lambda^{1/2}\left(\partial_{x}  \Delta_{\alpha} f\cos(\delta\Delta_{\alpha}f(x) ) \right) \ d\delta \ d\alpha \ dx  \\
&=& \int \Lambda f \ \int \partial_{x}  \Delta_{\alpha} f  \  \ \int_{0}^{\infty}   e^{-\delta}  \cos(\delta\Delta_{\alpha}f(x) ) \ d\delta \ d\alpha \ dx  \equiv L_1 
\end{eqnarray*}
 Let us set $\bar \Delta_{\alpha} f= \frac{f(x,t)-f(x+\alpha,t)}{\alpha}$. By denoting $S= \Delta_{\alpha} f + \bar \Delta_{\alpha} f$ and $D= \Delta_{\alpha} f - \bar \Delta_{\alpha} f,$ one observes that by doing $\alpha \rightarrow -\alpha$ if necessary, one may write
\begin{eqnarray*}
\frac{1}{2}  \partial_{t} \Vert f \Vert^{2}_{\dot H^{1/2}}&=&\int \Lambda f \ \int \partial_{x}D  \  \ \int_{0}^{\infty}  e^{-\delta}  \cos(\delta\Delta_{\alpha}f(x) ) \ d\delta \ d\alpha \ dx \\
&\quad &- \int \Lambda f \ \int \partial_{x}  \Delta_{\alpha} f \  \ \int_{0}^{\infty}   e^{-\delta}  \cos(\delta\bar\Delta_{\alpha}f ) \ d\delta \ d\alpha \ dx \\
&=&\int \Lambda f \ \int \partial_{x}D   \  \ \int_{0}^{\infty}  e^{-\delta}  \cos(\delta\Delta_{\alpha}f(x) ) \ d\delta \ d\alpha \ dx \\
&\quad &+ \int \Lambda f \ \int \partial_{x}  \Delta_{\alpha} f \  \ \int_{0}^{\infty}   e^{-\delta} \left( \cos(\delta\Delta_{\alpha}f) -\cos(\delta\bar\Delta_{\alpha}f) \right)\ d\delta \ d\alpha \ dx \\
&\quad &- \int \Lambda f \ \int \partial_{x}  \Delta_{\alpha} f \  \ \int_{0}^{\infty}   e^{-\delta}  \cos(\delta\Delta_{\alpha}f) \ d\delta \ d\alpha \ dx \\
&=&\frac{1}{4} \int \Lambda f  \int \partial_{x} D  \int_{0}^{\infty}  e^{-\delta} ( \cos(\delta\Delta_{\alpha}f(x))+\cos(\delta \bar\Delta_{\alpha}f(x)) \  d\delta \  d\alpha \  dx \\
&\quad &+ \frac{1}{2} \int \Lambda f  \int \partial_{x}  \Delta_{\alpha} f \int_{0}^{\infty}   e^{-\delta} \left(  \cos(\delta\Delta_{\alpha}f) -\cos(\delta\bar\Delta_{\alpha}f)  \right)\  d\delta \ d\alpha \ dx  \\
&=& \frac{1}{2}\int \Lambda f \int  \partial_{x}  D   \int_{0}^{\infty}   e^{-\delta} \cos(\frac{\delta}{2}D) \cos(\frac{\delta}{2}S) \ d\delta \ d\alpha \  dx \\
&\quad &- \int \Lambda f \ \int \partial_{x}  \Delta_{\alpha} f  \int_{0}^{\infty}   e^{-\delta}  
\sin(\frac{\delta}{2}D) \sin(\frac{\delta}{2}S)  \ d\delta \ d\alpha \ dx 
\end{eqnarray*}
Hence, we have that
\begin{eqnarray*}
\frac{1}{2}  \partial_{t} \Vert f \Vert^{2}_{\dot H^{1/2}}&=& -\int \Lambda f  \int \partial_{x}  D  \int_{0}^{\infty}   e^{-\delta} \cos(\frac{\delta}{2}D) \sin^{2}(\frac{\delta}{4}S)  \ d\delta \  d\alpha \  dx \\
&\quad&+\frac{1}{2} \int \Lambda f  \ \int \partial_{x}  D  \  \ \int_{0}^{\infty}   e^{-\delta} \cos(\frac{\delta}{2}D)\ d\delta \ d\alpha \ dx \\
&\quad&- \int \Lambda f \ \int \partial_{x}  \Delta_{\alpha} f \  \ \int_{0}^{\infty}   e^{-\delta}  
\sin(\frac{\delta}{2}D) \sin(\frac{\delta}{2}S)  \ d\delta \ d\alpha \ dx \\
&=& L_{1,1}+L_{1,2}+L_{1,3}
\end{eqnarray*}
We need to further decompose the last term, more precisely we write
\begin{eqnarray*}
\displaystyle L_{1,3}&=& - \int \Lambda f \ \int \partial_{x}  \Delta_{\alpha} f \  \ \int_{0}^{\infty}   e^{-\delta}  \sin(\frac{\delta}{2}D) \sin(\frac{\delta}{2} S)  \ d\delta \ d\alpha \ dx  \\
&=&- \int \Lambda f \ \int \frac{f_{x}(x)-f_{x}(x-\alpha)}{\alpha} \  \ \int_{0}^{\infty}   e^{-\delta}  \sin(\frac{\delta}{2}D) \sin(\frac{\delta}{2} S)  \ d\delta \ d\alpha \ dx \\
&=& \int \Lambda f \ \int \frac{f_{x}(x-\alpha)-f_{x}(x)}{\alpha} \  \ \int_{0}^{\infty}   e^{-\delta}  \sin(\frac{\delta}{2}D) \sin(\frac{\delta}{2} S)  \ d\delta \ d\alpha \ dx \\
&=& \int \Lambda f \ \int \frac{\partial_{\alpha}\left(f(x)-f(x-\alpha)\right)}{\alpha}  \int_{0}^{\infty}   e^{-\delta}  \sin(\frac{\delta}{2}D) \sin(\frac{\delta}{2} S)  \ d\delta \ d\alpha \ dx \\
&\quad&- \int \Lambda f \ \int \frac{f_{x}(x)}{\alpha} \int_{0}^{\infty}   e^{-\delta}  \sin(\frac{\delta}{2}D) \sin(\frac{\delta}{2} S)  \ d\delta \ d\alpha \ dx 
\end{eqnarray*}
We obtain, 
\begin{eqnarray*}
L_{1,3}&=& \int \Lambda f \int \frac{f(x)-f(x-\alpha)}{\alpha^2} \int_{0}^{\infty}   e^{-\delta}  
\sin(\frac{\delta}{2}D) \sin(\frac{\delta}{2}S)  \ d\delta \ d\alpha \ dx \\
&\quad&- \frac{1}{2}  \int \Lambda f  \int \frac{f(x)-f(x-\alpha)}{\alpha} \int_{0}^{\infty} \delta  e^{-\delta}  
\partial_{\alpha} D \cos(\frac{\delta}{2}D) \sin(\frac{\delta}{2}S)  \ d\delta \ d\alpha \ dx \\
&\quad&-\frac{1}{2}  \int \Lambda f  \int \frac{f(x)-f(x-\alpha)}{\alpha}  \int_{0}^{\infty} \delta  e^{-\delta}  
 \partial_{\alpha} S \sin(\frac{\delta}{2}D)  \cos(\frac{\delta}{2}S)  \ d\delta \ d\alpha \ dx \\
 &\quad&- \int \Lambda f \int \frac{f_{x}(x)}{\alpha} \int_{0}^{\infty}   e^{-\delta}  
\sin(\frac{\delta}{2}D) \sin(\frac{\delta}{2}S)  \ d\delta \ d\alpha \ dx  \\
&=& \sum_{i=1}^4 L_{1,3,i}
\end{eqnarray*}

\noindent{\bf{6.1.1. Estimates of $L_{1,1}$}} \\

\noindent In order to control $L_{1,1}$, we use Holder inequality $L^{2}-L^{2}-L^{\infty}$, 
\begin{eqnarray*}
L_{1,1}&=&-\int \Lambda f  \int \partial_{x}D    \int_{0}^{\infty}   e^{-\delta} \cos(\frac{\delta}{2}D) \sin^{2}(\frac{\delta}{4}S) \ d\delta \ d\alpha \ dx \\
&\leq&  \frac{\Gamma(3)}{4}\Vert f \Vert^{2}_{\dot H^{1}} \int \frac{\Vert \delta_{\alpha}f+\bar\delta_{\alpha}f \Vert^{2}_{L^{\infty}}}{\vert \alpha \vert^{3}} \ d\alpha \\
&\leq& \frac{1}{2} \Vert f \Vert^{2}_{\dot H^{1}} \Vert f \Vert^{2}_{\dot B^{1}_{\infty,2}} \\
&\leq& \frac{1}{2} \Vert f \Vert^{2}_{\dot H^{1}} \Vert f \Vert^{2}_{\dot H^{3/2}} \\
 \end{eqnarray*}
 where we used the embedding $\dot H^{3/2} \hookrightarrow \dot B^{1}_{\infty,2}$ along with the fact that (since $\dot H^{1}$ and $L^{2}$ are shift invariant spaces) $$\Vert \partial_{x} D \Vert_{L^{2}} \leq \frac{2}{\vert \alpha \vert}\Vert f_{x}(x)-f_{x}(x-\alpha)\Vert_{L^{2}} \leq \frac{4}{\vert \alpha \vert}\Vert \partial_{x} D \Vert_{L^{2}} \leq \frac{4}{\vert \alpha \vert}\Vert f\Vert_{\dot H^{1}}.$$ \\
 
 \noindent{\bf{6.1.2. Estimates of $L_{1,3}$}} \\

In the next subsection, we shall estimate the $L_{1,3,i}$ for $i=1,..,4$. \\

 \noindent{\bf{6.1.2.1. Estimates of $L_{1,3,1}$}} \\
 
By observing that $\dot H^{1} \hookrightarrow \dot B^{1/2}_{\infty,2}$, we have that
 \begin{eqnarray*}
\vert L_{1,3,1} \vert &\leq&  \Vert f \Vert_{\dot H^{1}}  \  \ \int_{0}^{\infty} \delta  e^{-\delta} \frac{\Vert f(x)-f(x-\alpha)\Vert_{L^{\infty}}\Vert f(x-\alpha)+f(x+\alpha)-2f(x)\Vert_{L^{2}}}{\vert \alpha \vert^{3}} \ d\delta \ d\alpha  \\
&\leq&  \Gamma(2) \Vert f \Vert_{\dot H^{1}} \left(\int \frac{\Vert f(x)-f(x-\alpha) \Vert^{2}_{L^{\infty}}}{\vert \alpha \vert^2} \ d\alpha \int \frac{\Vert f(x-\alpha)+f(x+\alpha)-2f(x) \Vert^{2}_{L^{2}}}{\vert \alpha \vert^4} \ d\alpha \right)^{1/2} \\
&\leq& \Vert f \Vert_{\dot H^{1}} \Vert f \Vert_{\dot H^{1/2}} \Vert f \Vert_{\dot B^{1/2}_{\infty,2}} \\
&\leq&  \Vert f \Vert^2_{\dot H^{1}} \Vert f \Vert_{\dot H^{3/2}}
 \end{eqnarray*}

  \noindent{\bf{6.1.2.2. Estimates of $L_{1,3,2}$}} \\
  
  We may rewrite (just by a direct integration in $x$) $D$ and $S$ as follows:  
  \begin{equation} \label{diff}
   D=\frac{f(x+\alpha)-f(x-\alpha)}{\alpha}=\frac{1}{\alpha} \int_0^\alpha (f_x(x+s)+f_x(x-s)-2f_x(x)) \ ds +2 f_x(x)
   \end{equation}
and
$S=\Delta_{\alpha}f+\bar\Delta_{\alpha}f=-\frac{(f(x+\alpha)+f(x-\alpha)-2f(x))}{\alpha},$ \\

We also need to give a suitable expression of their derivatives with respect to $\alpha$. In that regards, it is not difficult to check that
$$\partial_{\alpha}D=\frac{f_x(x+\alpha)+f_x(x-\alpha)-2f_x(x)}{\alpha}-\frac{\int_{0}^{\alpha} \left(f_x(x-s)+f_x(x+s)-2f_x(x)\right) \ ds}{\alpha^2}$$
and
$$\partial_{\alpha}S=\bar\Delta_{\alpha} f_x-\Delta_{\alpha} f_x+\frac{f(x+\alpha)+f(x-\alpha)-2f(x)}{\alpha^2},$$
we then rewrite $L_{1,3,2}$ as  
\begin{eqnarray*}
L_{1,3,2}&=& -\frac{1}{2} \int \Lambda f \ \int \frac{f(x)-f(x-\alpha)}{\alpha} \  \ \int_{0}^{\infty} \delta  e^{-\delta}  
\frac{f_x(x+\alpha)+f_x(x-\alpha)-2f_x(x)}{\alpha} \\
&& \times \ \cos(\frac{\delta}{2}(\Delta_{\alpha}f-\bar\Delta_{\alpha}f)) \sin(\frac{\delta}{2}(\Delta_{\alpha}f+\bar\Delta_{\alpha}f))  \ d\delta \ d\alpha \ dx \\
&\quad&+\frac{1}{2}  \int \Lambda f \ \int \frac{f(x)-f(x-\alpha)}{\alpha} \  \ \int_{0}^{\infty} \delta  e^{-\delta}  
 \ \frac{\int_{0}^{\alpha} \left(f_x(x-s)+f_x(x+s)-2f_x(x)\right) \ ds}{\alpha^2} \\
 &&\times \cos(\frac{\delta}{2}(\Delta_{\alpha}f-\bar\Delta_{\alpha}f)) \sin(\frac{\delta}{2}(\Delta_{\alpha}f+\bar\Delta_{\alpha}f))  \ d\delta \ d\alpha \ dx \\
 &=& L_{1,3,2,1} + L_{1,3,2,2}.
\end{eqnarray*}

By using that   $\dot H^1 \hookrightarrow \dot B^{1/2}_{\infty,2}$,  we find
\begin{eqnarray*}
\vert  L_{1,3,2,1} \vert &\leq& \frac{\Gamma(2)}{2}\Vert f \Vert_{\dot H^{1}} \int \frac{\Vert f(x)-f(x-\alpha) \Vert_{L^{\infty}}}{\alpha} \   
\frac{\Vert f_{x}(x+\alpha)+f_{x}(x-\alpha)-2f_{x}(x) \Vert_{L^{2}}}{\alpha} \ d\alpha  \\
&\leq& \Vert f \Vert_{\dot H^{1}} \Vert f \Vert_{\dot B^{1/2}_{\infty,2}} \Vert f_{x} \Vert_{\dot B^{1/2}_{2,2}} \\
&\leq& \Vert f \Vert^{2}_{\dot H^{1}} \Vert f\Vert_{\dot H^{3/2}}
\end{eqnarray*}
In order to estimate   $L_{1,3,2,2}$, we consider $q$, $r$ and $\bar r$ (so that $1/r+1/\bar r=1$) that will be chosen latter, and we write
\begin{eqnarray*}
\vert  L_{1,3,2,2} \vert &\leq&\frac{1}{2}\Vert f \Vert_{\dot H^{1}} \ \int \frac{\Vert f(x)-f(x-\alpha) \Vert_{L^{\infty}}}{\vert \alpha \vert^{3}} \  \ \int_{0}^{\infty} \delta  e^{-\delta} \\
&&\times   \vert \alpha \vert^{q+\frac{1}{\bar r}} \left(\int_{0}^{\alpha} \frac{ \Vert f_x(x-s)+f_x(x+s)-2f_x(x) \Vert^{r}_{L^{2}}}{s^{qr}} \ ds\right)^{1/r}\\
 &&\times  \frac{\Vert \delta_{\alpha}f+\bar\delta_{\alpha}f \Vert_{L^{\infty}}}{\alpha}  \ d\delta \ d\alpha  \\
 &\leq&\frac{\Gamma(2)}{2}\Vert f \Vert_{\dot H^{1}} \Vert f_x \Vert_{ \dot B^{q-\frac{1}{r}}_{2,r}} \ \int \frac{\Vert f(x)-f(x-\alpha) \Vert_{L^{\infty}}}{\vert \alpha \vert^{3-q-\frac{1}{\bar r}}} \  \frac{\Vert \delta_{\alpha}f+\bar\delta_{\alpha}f \Vert_{L^{\infty}}}{\alpha} \ d\alpha  \\
 &\leq& \frac{1}{2}\Vert f \Vert_{\dot H^{1}} \Vert f\Vert_{ \dot B^{1+q-\frac{1}{r}}_{2,r}} \left(\int \frac{\Vert f(x)-f(x-\alpha) \Vert^{2}_{L^{\infty}}}{\vert \alpha \vert^{5-2q-\frac{2}{\bar r}}} \ d\alpha \int \frac{\Vert \delta_{\alpha}f+\bar\delta_{\alpha}f \Vert^{2}_{L^{\infty}}}{\alpha^{3}} \ d\alpha \right)^{1/2} \\
 &\leq& \frac{1}{2}\Vert f \Vert_{\dot H^{1}} \Vert f \Vert_{ \dot B^{1+q-\frac{1}{r}}_{2,r}} \Vert f \Vert_{\dot B^{2-q-\frac{1}{\bar r}}_{\infty,2}}\Vert f \Vert_{\dot B^{1}_{\infty,2}} 
 \end{eqnarray*}
Then, for $\bar r= r=2$ and $q=1$, we obtain since $\dot H^{1} \hookrightarrow \dot B^{1/2}_{\infty,2}$ and $\dot H^{3/2} \hookrightarrow \dot B^{1}_{\infty,2}$
\begin{eqnarray*}
\vert  L_{1,3,2,2} \vert &\leq& \frac{1}{2}\Vert f \Vert_{\dot H^{1}} \Vert f \Vert_{ \dot H^{3/2}} \Vert f \Vert_{\dot B^{1/2}_{\infty,2}} \Vert f \Vert_{\dot B^{1}_{\infty,2}} 
 \\
  &\leq&  \frac{1}{2}\Vert f \Vert^{2}_{\dot H^{1}} \Vert f \Vert^{2}_{\dot H^{3/2}},
\end{eqnarray*}
hence,
\begin{eqnarray*} \vert L_{1,3,2} \vert \leq \Vert f \Vert^{2}_{\dot H^{1}} \Vert f \Vert^{2}_{\dot H^{3/2}}
\end{eqnarray*}

  \noindent{\bf{6.1.2.3. Estimates of $L_{1,3,3}$}} \\

\noindent We now estimate $L_{1,3,3}$. We need to decompose this term as follows

\begin{eqnarray*}
  L_{1,3,3} &=&-\frac{1}{2}\int \Lambda f  \int \frac{f(x)-f(x-\alpha)}{\alpha} \int_{0}^{\infty} \delta  e^{-\delta}   \sin(\frac{\delta}{2} D_{\alpha}f) \partial_{\alpha} S  \cos(\frac{\delta}{2}S)  \  d\delta \ d\alpha \ dx \\
 &=&  -\frac{1}{2}\int \Lambda f  \int \frac{f(x)-f(x-\alpha)}{\alpha} \ \int_{0}^{\infty} \delta  e^{-\delta}  
 \ \sin(\frac{\delta}{2}D) \bar \Delta_{\alpha} f_{x} \ \cos(\frac{\delta}{2}S)   \  d\delta \ d\alpha \ dx \\
&\quad&+\frac{1}{2} \int \Lambda f  \int \frac{f(x)-f(x-\alpha)}{\alpha}   \ \int_{0}^{\infty} \delta  e^{-\delta} 
 \ \sin(\frac{\delta}{2}D) \Delta_{\alpha}f_{x} \ \cos(\frac{\delta}{2}S)  \  d\delta \ d\alpha \ dx \\ &\quad&-\frac{1}{2} \int \Lambda f \ \int \frac{f(x)-f(x-\alpha)}{\alpha}  \int_{0}^{\infty} \delta  e^{-\delta}   \sin(\frac{\delta}{2}D) \frac{f(x+\alpha)+f(x-\alpha)-2f(x)}{\alpha^2} \\
 &\quad& \times\cos(\frac{\delta}{2}S)  \ d\delta \ d\alpha \ dx \\
  &=& \sum_{i=1}^{3}    L_{1,3,3, i} 
  \end{eqnarray*}
    
  \noindent To estimate $L_{1,3,3,1}$ we see that 
  
  \begin{eqnarray*}
  \vert L_{1,3,3,1} \vert &\leq&\frac{\Gamma(2)}{2} \Vert f \Vert_{\dot H^{1}} \left(\int \frac{\Vert f(x)-f(x-\alpha) \Vert^{2}_{L^{\infty}}}{\vert \alpha \vert^{2}} \ d\alpha \int \frac{\Vert f_{x}(x)-f_{x}(x+\alpha) \Vert^{2}_{L^{2}}}{\vert \alpha \vert^{2}} \ d\alpha \right)^{1/2} \\
   &\leq&  \Vert f \Vert_{\dot H^{1}} \Vert f \Vert_{\dot B^{1/2}_{\infty,2}} \Vert f_{x} \Vert_{\dot B^{1/2}_{2,2}} \\
    &\leq& \Vert f \Vert^{2}_{\dot H^{1}} \Vert f \Vert_{\dot H^{3/2}}
  \end{eqnarray*}
  As well, one may easily estimate $L_{1,3,3,2}$ and we find that
   \begin{eqnarray*}
   \vert L_{1,3,3,2} \vert  &\leq& \frac{1}{2}\Vert f \Vert^{2}_{\dot H^{1}} \Vert f \Vert_{\dot H^{3/2}}
  \end{eqnarray*}
  For $L_{1,3,3,3}$, it suffices to write that
  \begin{eqnarray*}
   \vert L_{1,3,3,3} \vert &\leq& \frac{1}{2}\Vert f \Vert_{\dot H^{1}} \left(\int \frac{\Vert f(x)-f(x-\alpha) \Vert^{2}_{L^{\infty}}}{\vert \alpha \vert^{2}} \ d\alpha \int \frac{\Vert f(x+\alpha)+f(x-\alpha)-2f(x)\Vert^{2}_{L^{2}}}{\vert \alpha \vert^{4}} \ d\alpha \right)^{1/2} \\
   &\leq&\frac{1}{2}\Vert f \Vert_{\dot H^{1}} \Vert f \Vert_{\dot B^{1/2}_{\infty,2}}\Vert f \Vert_{\dot B^{3/2}_{2,2}} \\
      &\leq&\frac{1}{2}\Vert f \Vert^{2}_{\dot H^{1}} \Vert f \Vert_{\dot H^{3/2}}
   \end{eqnarray*}
   So that,
   \begin{equation} \label{eq:croisŽs2}
    \vert L_{1,3,3} \vert \lesssim \Vert f \Vert^{2}_{\dot H^{1}} \Vert f \Vert_{\dot H^{3/2}}
   \end{equation}
   
       \noindent{\bf{6.1.2.4. Estimates of $L_{1,3,4}$}} \\

As for $L_{1,3,4}$, we use the antisymmetry of $\mathcal{H}$ together with the commutator estimates \eqref{cz}, as follows
\begin{eqnarray*}
L_{1,3,4} &=& - \int \Lambda f \ \int \frac{f_{x}(x)}{\alpha} \  \ \int_{0}^{\infty}   e^{-\delta}  
\sin(\frac{\delta}{2}D) \sin(\frac{\delta}{2}S)  \ d\delta \ d\alpha \ dx \\
&=& \frac{1}{2} \int f_{x} \ \int_{0}^{\infty} \int    e^{-\delta}  \frac{1}{\alpha}
 \left[\mathcal{H} , \sin(\frac{\delta}{2}D) \sin(\frac{\delta}{2}S) \right]f_{x}   \ d\delta \ d\alpha \ dx. 
     \end{eqnarray*}
and we find that
 \begin{eqnarray*}
 L_{1,3,4} &=&  \frac{1}{2}\int\int \frac{f_{x}(x)-f_{x}(x-\alpha)}{\alpha} \int_{0}^{\infty}     e^{-\delta}  
 \left[\mathcal{H} , \sin(\frac{\delta}{2}D) \sin(\frac{\delta}{2}S) \right]f_{x}  \ d\delta \ d\alpha \ dx \\
  &\quad&+\frac{1}{2} \int \int \frac{f_{x}(x-\alpha)}{\alpha} \ \int_{0}^{\infty}     e^{-\delta}  
 \left[\mathcal{H} , \sin(\frac{\delta}{2}D) \sin(\frac{\delta}{2}S) \right]f_{x}   \ d\delta \ d\alpha \ dx.
   \end{eqnarray*}
   By integrating by parts we find
   \begin{eqnarray*}
L_{1,3,4} &=& \frac{1}{2}  \int \int \frac{f_{x}(x)-f_{x}(x-\alpha)}{\alpha}  \int_{0}^{\infty}     e^{-\delta}  
\left[\mathcal{H} , \sin(\frac{\delta}{2}D) \sin(\frac{\delta}{2}S) \right]f_{x} \  d\delta \ d\alpha \ dx   \\ 
&\quad&-\frac{1}{2} \int \int \frac{f(x-\alpha)-f(x)}{\alpha^2} \ \int_{0}^{\infty}     e^{-\delta}  
 \left[\mathcal{H} , \sin(\frac{\delta}{2}D) \sin(\frac{\delta}{2}S) \right]f_{x} \ d\delta \ d\alpha \ dx   \\
 &\quad&+ \frac{1}{2}\int\int \frac{f(x-\alpha)-f(x)}{\alpha} \ \int_{0}^{\infty}   e^{-\delta}  
 \partial_{\alpha}\left[\mathcal{H} , \sin(\frac{\delta}{2}D) \sin(\frac{\delta}{2}S) \right]f_{x}  \  d\delta \ d\alpha \ dx   \\
 &=& L_{1,3,4,1}+L_{1,3,4,2}+L_{1,3,4,3}
   \end{eqnarray*}
  The commutator estimate \eqref{cz} (in the case $l=1$ and $k=0$) and the embeddings $\dot H^{3/2} \hookrightarrow \dot B^{1}_{\infty,4}$ and $\dot H^{1} \hookrightarrow \dot B^{1/2}_{\infty,4}$ allows us to find
\begin{eqnarray*}
 \vert L_{1,3,4,1} \vert &\leq& \frac{\Gamma(2)}{2} \Vert f\Vert_{\dot H^{1}} \int  \frac{\Vert f_{x}(x)-f_{x}(x-\alpha) \Vert_{L^{2}}}{\alpha} \frac{\Vert  f(x)-f(x-\alpha)\Vert_{L^{\infty}}}{\alpha}  \\
 &&\hspace{2cm} \times\ \frac{\Vert f (x-\alpha)+f(x+\alpha)-2f(x) \Vert_{L^{\infty}}}{\alpha}  \ d\alpha \\
 &\leq&  \Vert f\Vert_{\dot H^{1}}   \left(\int \frac{\Vert f_{x}(x)-f_{x}(x-\alpha) \Vert^{2}_{L^{2}}}{\alpha^2} \ d\alpha \right)^{1/2} \left(\int \frac{\Vert f (x-\alpha)+f(x+\alpha)-2f(x) \Vert^{4}_{L^{\infty}}}{\alpha^{5}}  \ d\alpha \right)^{1/4} \\
 &\quad&\times  \left(\int\frac{\Vert  f(x)-f(x-\alpha)\Vert^{4}_{L^{\infty}}}{\alpha^3} \ d\alpha \right)^{1/4}  \\
 &\leq&  \Vert f\Vert_{\dot H^{1}} \Vert f \Vert_{\dot H^{3/2}} \Vert f \Vert_{\dot B^{1}_{\infty,4}} \Vert f \Vert_{\dot B^{1/2}_{\infty,4}} \\
 &\leq&  \Vert f\Vert^{2}_{\dot H^{1}} \Vert f \Vert^{2}_{\dot H^{3/2}} 
 \end{eqnarray*}
Then, we estimate $L_{1,3,4,2}$. We first see that
\begin{eqnarray*}
L_{1,3,4,2}&=& -\frac{1}{2}\int \Lambda f \frac{f(x-\alpha)-f(x)}{\alpha^2}  \int_{0}^{\infty} \int    e^{-\delta}  \sin(\frac{\delta}{2}D) \sin(\frac{\delta}{2}S) \  d\delta  \ d\alpha \ dx \\
&\quad&-\frac{1}{2} \int f_{x} \frac{\mathcal{H}f(x-\alpha)-\mathcal{H}f(x)}{\alpha^2}  \int_{0}^{\infty} \int    e^{-\delta}  \sin(\frac{\delta}{2}D) \sin(\frac{\delta}{2}S)  \  d\delta \ d\alpha \ dx \\
&=& L_{1,3,4,2,1}+ L_{1,3,4,2,2}
\end{eqnarray*}
Then, we have  by using Holder inequality $L^{2}-L^{\infty}-L^{2}$
\begin{eqnarray*}
\vert L_{1,3,4,2,1} \vert &\leq&\frac{1}{2} \Vert f \Vert_{\dot H^{1}} \int \frac{\Vert f(x-\alpha)-f(x) \Vert_{L^{\infty}}}{\alpha^2} \ \int_{0}^{\infty}     e^{-\delta}  \frac{\Vert \delta_{\alpha}f+\bar\delta_{\alpha}f \Vert_{L^{2}}}{\alpha}  \  \ d\delta \ d\alpha   \\
&\leq& \frac{1}{2} \Vert f \Vert_{\dot H^{1}} \left(\int \frac{\Vert f(x-\alpha)-f(x) \Vert^{2}_{L^{\infty}}}{\alpha^3} \ d\alpha \right)^{1/2}  \left( \int \frac{\Vert \delta_{\alpha}f+\bar\delta_{\alpha}f \Vert^{2}_{L^{2}}}{\alpha^3} \ d\alpha \right)^{1/2} \\
&\leq&\frac{1}{2} \Vert f \Vert_{\dot H^{1}}  \Vert f \Vert_{\dot B^{1}_{\infty,2}} \Vert f \Vert_{\dot H^{1}} \\
&\leq&\frac{1}{2} \Vert f \Vert^{2}_{\dot H^{1}} \Vert f \Vert_{\dot H^{3/2}}  
\end{eqnarray*}
Analogously  (since $\mathcal{H}$ maps $L^{2}$ onto $L^{2}$), we also find
\begin{equation*}
\vert L_{1,3,4,2} \vert \leq \frac{1}{2} \Vert f \Vert^{2}_{\dot H^{1}} \Vert f \Vert_{\dot H^{3/2}}
\end{equation*}
The control of the term $L_{1,3,4,3}$ is challenging and one needs to use the following decomposition
\begin{eqnarray*}
 L_{1,3,4,3}&=&-\frac{1}{4}\int f_{x} \frac{\mathcal{H}f(x-\alpha)-\mathcal{H}f(x)}{\alpha} \ \int_{0}^{\infty} \int  \delta  e^{-\delta}  \partial_{\alpha}D \ \cos(\frac{\delta}{2}D) \sin(\frac{\delta}{2}S)   \ d\delta \ d\alpha \ dx \\
 &\quad&-\frac{1}{4}\int \Lambda f \frac{f(x-\alpha)-f(x)}{\alpha} \ \int_{0}^{\infty} \int  \delta  e^{-\delta}  \partial_{\alpha}D \ \cos(\frac{\delta}{2}D) \sin(\frac{\delta}{2}S)    \ d\delta \ d\alpha \ dx \\
 &\quad&-\frac{1}{4} \int f_{x} \int \frac{\mathcal{H}f(x-\alpha)-\mathcal{H}f(x)}{\alpha} \ \int_{0}^{\infty}   \delta  e^{-\delta}  \partial_{\alpha}S \ \cos(\frac{\delta}{2}S)  \sin(\frac{\delta}{2}D)  \  \ d\delta \ d\alpha \ dx \\
 &\quad&-\frac{1}{4} \int \Lambda f \int \frac{f(x-\alpha)-f(x)}{\alpha} \ \int_{0}^{\infty}   \delta  e^{-\delta}  \partial_{\alpha}S \ \cos(\frac{\delta}{2}S)  \sin(\frac{\delta}{2}D)  \  \ d\delta \ d\alpha \ dx \\
 &=& \sum_{i=1}^{4}  L_{1,3,4,3,i}
\end{eqnarray*}
All those previous term are as regular as $L_{1,3,j}$ for $j=2,3$ (up to some Hilbert transform, we shall do $L^{p}$ estimate for safe values of $p$, that is $p\ne 1, \infty$), it is therefore an easy task to get that 
$$
\left \vert L_{1,3,4,3,1}+L_{1,3,4,3,2} \right\vert \leq 2 \Vert f \Vert^{2}_{\dot H^{1}} \Vert f \Vert^{2}_{\dot H^{3/2}}
$$
and
$$
\left \vert L_{1,3,4,3,3}+L_{1,3,4,3,4} \right\vert \leq 3\Vert f \Vert^{2}_{\dot H^{1}} \Vert f \Vert_{\dot H^{3/2}} \\
$$

Therefore, we have obtained

\begin{equation} \label{i13}
\vert L_{1,3} \vert \leq \Vert f \Vert^{2}_{\dot H^{1}} \left(\Vert f \Vert^{2}_{\dot H^{3/2}} +  \Vert f \Vert_{\dot H^{3/2}} \right).
\end{equation}

\noindent {\bf{6.1.3. Estimates of $L_{1,2}$}}  \\
 
\noindent We have
 \begin{eqnarray*}
 L_{1,2}&=&\frac{1}{2}\int \Lambda f \ \int (\partial_{x}  \Delta_{\alpha} f -\partial_{x} \bar \Delta_{\alpha} f)  \  \ \int_{0}^{\infty}   e^{-\delta} \cos(\frac{\delta}{2}D)\ d\delta \ d\alpha \ dx \\
  &=& -\frac{1}{2}\int \Lambda f \ \int \frac{1}{\alpha} \partial_{\alpha}(\delta_{\alpha} f + \bar \delta_{\alpha} f )  \  \ \int_{0}^{\infty}   e^{-\delta} \cos(\frac{\delta}{2}D)\ d\delta \ d\alpha \ dx, \\
  \end{eqnarray*}
By integrating by parts one finds
   \begin{eqnarray*}
L_{1,2}  &=&\frac{1}{2}  \int \Lambda f  \int \frac{f(x-\alpha)+f(x+\alpha)-2f(x)}{\alpha^2}    \  \ \int_{0}^{\infty}   e^{-\delta} \cos(\frac{\delta}{2} D)\ d\delta \ d\alpha \ dx \\
&\quad&+\frac{1}{2} \int \Lambda f  \int \frac{f(x-\alpha)+f(x+\alpha)-2f(x)}{\alpha}    \  \ \int_{0}^{\infty} \delta  e^{-\delta} \partial_{\alpha}D \sin(\frac{\delta}{2} D)  d\delta \ d\alpha \ dx. \\
  \end{eqnarray*}
 Therefore,
  \begin{eqnarray*}
L_{1,2} &=& \frac{1}{2}\int \Lambda f \ \int \frac{f(x-\alpha)+f(x+\alpha)-2f(x)}{\alpha^2}    \  \ \int_{0}^{\infty}   e^{-\delta} \cos(\frac{\delta}{2} D)\ d\delta \ d\alpha \ dx \\
&\quad&+\frac{1}{2} \int \Lambda f \ \int \frac{f(x-\alpha)+f(x+\alpha)-2f(x)}{\alpha}    \  \ \int_{0}^{\infty} \delta  e^{-\delta} \frac{f_x(x-\alpha)+f_x(x+\alpha)-2f_x(x)}{\alpha} \\
&\quad&\times\sin(\frac{\delta}{2} (\Delta_{\alpha}f-\bar\Delta_{\alpha}f))\ d\delta \ d\alpha \ dx \\
&\quad&+\frac{1}{2}\int \Lambda f\ \int \frac{f(x-\alpha)+f(x+\alpha)-2f(x)}{\alpha}    \  \ \int_{0}^{\infty} \delta  e^{-\delta} \frac{f_{x}(x)}{\alpha} \sin(\frac{\delta}{2}D)\ d\delta \ d\alpha \ dx \\
&\quad&-\frac{1}{2}\int \Lambda f \ \int \frac{f_{x}(x-\alpha)+f_{x}(x+\alpha)-2f_{x}(x)}{\alpha}    \  \ \int_{0}^{\infty} \delta  e^{-\delta} \frac{f(x+\alpha)-f(x-\alpha)}{\alpha^2} \\
&\quad&\times\sin(\frac{\delta}{2} (\Delta_{\alpha}f-\bar\Delta_{\alpha}f))\ d\delta \ d\alpha \ dx \\
\end{eqnarray*}

The first term gives the paraboliticity, since we may write \\

\begin{eqnarray*}
&&\frac{1}{2}\int \Lambda f \ \int \frac{f(x-\alpha)+f(x+\alpha)-2f(x)}{\alpha^2}    \int_{0}^{\infty}   e^{-\delta} \cos(\frac{\delta}{2} D)\ d\delta \ d\alpha \ dx \\
&&= -  \pi  \int \vert \Lambda f \vert^{2} \ dx - \int \Lambda f  \int \frac{f(x-\alpha)+f(x+\alpha)-2f(x)}{\alpha^2}    \int_{0}^{\infty}   e^{-\delta} \sin^{2}(\frac{\delta}{4} D)\ d\delta \ d\alpha \ dx, \\
&&= -  \pi \Vert f \Vert^{2}_{\dot H^{1}} - \int \Lambda f  \int \frac{f(x-\alpha)+f(x+\alpha)-2f(x)}{\alpha^2}    \int_{0}^{\infty}   e^{-\delta} \sin^{2}(\frac{\delta}{4} D)\ d\delta \ d\alpha \ dx,
\end{eqnarray*}

Then,  using formula \eqref{diff} we may rewrite this term as follows
\begin{eqnarray} \label{10sip}
L_{1,2}&=& \nonumber - \int \Lambda f \ \int \frac{f(x-\alpha)+f(x+\alpha)-2f(x)}{\alpha^2}    \  \ \int_{0}^{\infty}  e^{-\delta} \sin^{2}(\frac{\delta}{4} D)\ d\delta \ d\alpha \ dx \\
&\quad&+\nonumber \frac{1}{2} \int \Lambda f \ \int \frac{f(x-\alpha)+f(x+\alpha)-2f(x)}{\alpha}    \  \ \int_{0}^{\infty} \delta  e^{-\delta} \frac{f_{x}(x-\alpha)+f_{x}(x+\alpha)-2f_{x}(x)}{\alpha} \\
&\quad&\times \nonumber  \sin(\frac{\delta}{2} (\Delta_{\alpha}f-\bar\Delta_{\alpha}f))\ d\delta \ d\alpha \ dx \\
&\quad&-\nonumber \frac{1}{2}\int \Lambda f \ \int \frac{f(x-\alpha)+f(x+\alpha)-2f(x)}{\alpha^3}    \  \ \int_{0}^{\infty} \delta  e^{-\delta} \int_{0}^{\alpha} f_{x}(x-s)+f_{x}(x+s)-2f_{x}(x) \ ds \\
&\quad&\times\nonumber \sin(\frac{\delta}{2} (\Delta_{\alpha}f-\bar\Delta_{\alpha}f))\ d\delta \ d\alpha \ dx \\
&\quad&-   \pi  \Vert f \Vert^{2}_{\dot H^{1}} \\
&=&\nonumber L_{1,2,1}+L_{1,2,2}+L_{1,2,3}+L_{1,2,4}.
 \end{eqnarray} 

We need to futher decompose $L_{1,2,1}$ as follows 
 
\begin{eqnarray*}
L_{1,2,1}&=&- \int \Lambda f \ \int \frac{f(x-\alpha)+f(x+\alpha)-2f(x)}{\alpha^2}    \  \ \int_{0}^{\infty}   e^{-\delta} 
\sin^{2}(\frac{\delta}{4} D)\ d\delta \ d\alpha \ dx \\
&=&- \int \Lambda f \ \int \frac{f(x-\alpha)+f(x+\alpha)-2f(x)}{\alpha^2}    \  \ \int_{0}^{\infty}   e^{-\delta} \\
&\quad&\times \sin(\frac{\delta}{4} \frac{1}{\alpha} \int_{0}^{\alpha} f_{x}(x+s)+f_{x}(x-s)-2f_{x}(x) \ ds - \frac{\delta}{2} f_{x}(x) )\ \sin(\frac{\delta}{4} D) \ d\delta \ d\alpha \ dx \\
&=&- \int \Lambda f \ \int \frac{f(x-\alpha)+f(x+\alpha)-2f(x)}{\alpha^2}    \  \ \int_{0}^{\infty}   e^{-\delta} \\
&\quad&\times \sin(\frac{\delta}{4} \frac{1}{\alpha} \int_{0}^{\alpha} f_{x}(x+s)+f_{x}(x-s)-2f_{x}(x) \ ds)\cos(\frac{\delta}{2} f_{x}(x))\ \sin(\frac{\delta}{4} D) \ d\delta \ d\alpha \ dx \\
&\quad&- \int \mathcal{H}f_{x} \ \int \frac{f(x-\alpha)+f(x+\alpha)-2f(x)}{\alpha^2}    \  \ \int_{0}^{\infty}  e^{-\delta} \\
&\quad&\times \cos\left(\frac{\delta}{4} \frac{1}{\alpha} \int_{0}^{\alpha} f_{x}(x+s)+f_{x}(x-s)-2f_{x}(x) \ ds\right)\sin(\frac{\delta}{2} f_{x}(x))\ \sin(\frac{\delta}{4} D) \ d\delta \ d\alpha \ dx 
\end{eqnarray*} 
In the last integral, with add and substract 1 in the cosine term, and we shall repeat this process after having used the trigonometry formula we develop the $\sin(\frac{\delta}{4} D)$, we obtain that
 \begin{eqnarray*}
L_{1,2,1}&=&-  \int \Lambda f \ \int \frac{f(x-\alpha)+f(x+\alpha)-2f(x)}{\alpha^2}    \  \ \int_{0}^{\infty}   e^{-\delta} \\
&\quad&\times \sin(\frac{\delta}{4} \frac{1}{\alpha} \int_{0}^{\alpha} f_{x}(x+s)+f_{x}(x-s)-2f_{x}(x) \ ds)\cos(\frac{\delta}{2} f_{x}(x))\ \sin(\frac{\delta}{4} D) \ d\delta \ d\alpha \ dx \\
&\quad&+2 \int \Lambda f \ \int \frac{f(x-\alpha)+f(x+\alpha)-2f(x)}{\alpha^2}    \  \ \int_{0}^{\infty}  e^{-\delta} \\
&\quad&\times \sin^{2} \left(\frac{\delta}{8} \frac{1}{\alpha} \int_{0}^{\alpha} f_{x}(x+s)+f_{x}(x-s)-2f_{x}(x) \ ds\right)\sin(\frac{\delta}{2} f_{x}(x))\ \sin(\frac{\delta}{4} D) \ d\delta \ d\alpha \ dx \\
&\quad&-\int \Lambda f \ \int \frac{f(x-\alpha)+f(x+\alpha)-2f(x)}{\alpha^2}    \  \ \int_{0}^{\infty}   e^{-\delta} \\
&\quad&\times \sin(\frac{\delta}{2} f_{x}(x))\ \sin(\frac{\delta}{4} \frac{1}{\alpha} \int_{0}^{\alpha} f_{x}(x+s)+f_{x}(x-s)-2f_{x}(x) \ ds)\cos( \frac{\delta}{2} f_{x}(x)) \ d\delta \ d\alpha \ dx  \\
&\quad&-2 \int \Lambda f \ \int \frac{f(x-\alpha)+f(x+\alpha)-2f(x)}{\alpha^2}    \  \ \int_{0}^{\infty}  e^{-\delta} \\
&\quad&\times \sin^{2}(\frac{\delta}{2} f_{x}(x))\ \sin^{2}\left(\frac{\delta}{8} \frac{1}{\alpha} \int_{0}^{\alpha} f_{x}(x+s)+f_{x}(x-s)-2f_{x}(x) \ ds\right) \ d\delta \ d\alpha \ dx \\
&\quad&- \int \Lambda f \ \int \frac{f(x-\alpha)+f(x+\alpha)-2f(x)}{\alpha^2}    \  \ \int_{0}^{\infty}   e^{-\delta} \ \sin^{2}( \frac{\delta}{2} f_{x}(x)) \ d\delta \ d\alpha \ dx \\
&=& \sum_{j=1}^{5} L_{1,2,1,j}
 \end{eqnarray*}
For $L_{1,2,1,1}$ we write 
 \begin{eqnarray*}
 \vert L_{1,2,1,1} \vert &\leq&  \Vert f \Vert_{\dot H^{1}} \int_{0}^{\infty} \delta e^{-\delta} \frac{\Vert f(x-\alpha)+f(x+\alpha)-2f(x) \Vert_{L^{\infty}}}{\vert \alpha \vert^{3}} \\
 && \ \ \times \int_{0}^{\alpha}  {\Vert f_{x}(x-s)+f_{x}(x+s)-2f_{x}(x) \Vert_{L^{2}}}  \ ds \ d\alpha \ d\delta 
 \end{eqnarray*}
 By  Minkowski's inequality, we find
 \begin{eqnarray*}
  \vert L_{1,2,1,1} \vert &\leq&  \Vert f \Vert_{\dot H^{1}} \int_{0}^{\infty} \delta e^{-\delta} \int \frac{\Vert f(x-\alpha)-f(x) \Vert_{L^{\infty}}}{\vert \alpha \vert^{3}}  \int_{0}^{\alpha}  {\Vert f_{x}(x-s)-f_{x}(x) \Vert_{L^{2}}}  \ ds \ d\alpha \ d\delta \\
 &\quad&+  \Vert f \Vert_{\dot H^{1}} \int_{0}^{\infty} \delta e^{-\delta} \int\frac{\Vert f(x-\alpha)-f(x) \Vert_{L^{\infty}}}{\vert \alpha \vert^{3}} \int_{0}^{\alpha}  {\Vert f_{x}(x+s)-f_{x}(x) \Vert_{L^{2}}}  \ ds \ d\alpha \ d\delta \\
 &\quad&+   \Vert f \Vert_{\dot H^{1}} \int_{0}^{\infty} \delta e^{-\delta} \int \frac{\Vert f(x+\alpha)-f(x) \Vert_{L^{\infty}}}{\vert \alpha \vert^{3}}  \int_{0}^{\alpha}  {\Vert f_{x}(x-s)-f_{x}(x) \Vert_{L^{2}}}  \ ds \ d\alpha \ d\delta \\
 &\quad&+   \Vert f \Vert_{\dot H^{1}} \int_{0}^{\infty} \delta e^{-\delta} \int\frac{\Vert f(x+\alpha)-f(x) \Vert_{L^{\infty}}}{\vert \alpha \vert^{3}}  \int_{0}^{\alpha}  {\Vert f_{x}(x+s)-f_{x}(x) \Vert_{L^{2}}}  \ ds \ d\alpha \ d\delta
 \end{eqnarray*}
 Since those terms have the same regularity, it is easy to conclude that
   \begin{eqnarray*}
 \vert L_{1,2,1,1} \vert &\leq& 4 \Gamma(2) \Vert f \Vert_{\dot H^{1}} \int \frac{\Vert f(x-\alpha)-f(x) \Vert_{L^{\infty}}}{\vert \alpha \vert^{3}}  \int_{0}^{\alpha}  {\Vert f_{x}(x-s)-f_{x}(x) \Vert_{L^{2}}} \ ds \ d\alpha  \\
 &\leq&  4 \Vert f \Vert_{\dot H^{1}} \int \frac{\Vert f(x-\alpha)-f(x) \Vert_{L^{\infty}}}{\vert \alpha \vert^{3}} \vert \alpha \vert^{q+\frac{1}{\bar r}}   \left(\int  \frac{\Vert f_{x}(x-s)-f_{x}(x) \Vert^{r} _{L^{2}}}{\vert s \vert^{qr}} \ ds\right)^{1/r}  \ d\alpha  \\
 &\leq&4 \Vert f \Vert_{\dot H^{1}} \Vert f \Vert_{\dot B^{2-q-\frac{1}{\bar r}}_{\infty,1}} \Vert f_{x} \Vert_{\dot B^{q-\frac{1}{r}}_{2,r}} 
  \end{eqnarray*}
 Then, choosing $q=3/4,$  $r=\bar r=2$ , using interpolation and the embedding $\dot H^{3/2} \hookrightarrow \dot B^{1}_{\infty,\infty}$,
  \begin{eqnarray*}
\vert L_{1,2,1,1} \vert  &\leq& 4 \Vert f \Vert_{\dot H^{1}} \Vert f \Vert_{\dot B^{\frac{3}{4}}_{\infty,1}} \Vert f_{x} \Vert_{\dot B^{1/4}_{2,2}} \\
&\leq& 4\Vert f \Vert_{\dot H^{1}} \Vert f \Vert^{1/2}_{\dot B^{1/2}_{\infty,\infty}} \Vert f \Vert^{1/2}_{\dot B^{1}_{\infty,\infty}} \Vert f \Vert_{\dot H^{5/4}} \\
&\leq&4 \Vert f \Vert_{\dot H^{1}} \Vert f \Vert^{1/2}_{\dot H^{1}} \Vert f \Vert^{1/2}_{\dot B^{1}_{\infty,\infty}} \Vert f \Vert^{1/2}_{\dot B^{1}_{2,2}}\Vert f \Vert^{1/2}_{\dot B^{3/2}_{2,2}} \\
 &\leq& 4 \Vert f\Vert^{2}_{\dot H^{1}} \Vert f\Vert_{\dot H^{3/2}}
  \end{eqnarray*} 
   Thus,
$$
\vert L_{1,2,1,1} \vert \leq\Vert f \Vert^{2}_{\dot H^{1}}  \Vert f \Vert_{\dot H^{3/2}}\\ 
$$

 \noindent Analogously, we find that
 
\begin{eqnarray}
\left\vert  \sum_{i=1}^{4} L_{1,2,1,i} \right \vert \leq \Vert f \Vert^{2}_{\dot H^{1}}  \Vert f \Vert_{\dot H^{3/2}} 
\end{eqnarray}

As for $L_{1,2,1,5}$, we observe that since
$$
-2 \pi \Lambda f = \int \frac{f(x-\alpha)+f(x+\alpha)-2f(x)}{\alpha^2} \ d\alpha \ \ \ {\rm{and}} \ \ \ \int_{0}^{\infty}   e^{-\delta} \ \sin^{2}(\frac{\delta}{2} A) \ d\delta = \frac{1}{2} \frac{1}{1+A^2}
$$

\begin{eqnarray} \label{a10sip}
 L_{1,2,1,5} &=&\nonumber -  \int \Lambda f \ \int \frac{f(x-\alpha)+f(x+\alpha)-2f(x)}{\alpha^2}    \  \ \int_{0}^{\infty}   e^{-\delta} \ \sin^{2}( \frac{\delta}{2} f_{x}(x)) \ d\delta \ d\alpha \ dx   \\
&\leq& \pi \frac{{K}^2}{1+K^2} \Vert f \Vert^{2}_{\dot H^1},
\end{eqnarray}

where 
$$
K=\Vert f_{x} \Vert_{L^{\infty}L^{\infty}}.
$$

Therefore, we find that 

\begin{equation}
\left\vert  \sum_{i=1}^{5} L_{1,2,1,i} \right \vert \leq \Vert f \Vert^{2}_{\dot H^{1}} \Vert f \Vert_{\dot H^{3/2}} +   \pi\frac{K^2}{1+K^2}  \Vert f \Vert^2_{\dot H^1}
\end{equation}

For  $L_{1,2,2}$, we write that

\begin{eqnarray*}
L_{1,2,2} &=& \frac{1}{2}\int \Lambda f \ \int \frac{f(x-\alpha)+f(x+\alpha)-2f(x)}{\alpha}    \  \ \int_{0}^{\infty} \delta  e^{-\delta} \frac{f_{x}(x-\alpha)+f_{x}(x+\alpha)-2f_{x}(x)}{\alpha} \\
&\quad&\times \sin(\frac{\delta}{2} (\Delta_{\alpha}f-\bar\Delta_{\alpha}f))\ d\delta \ d\alpha \ dx \\
&\leq&\frac{\Gamma(2)}{2} \Vert f \Vert_{\dot H^{1}}\left( \int 
\frac{\Vert f(x-\alpha)+f(x+\alpha)-2f(x)\Vert^{2}_{L^{\infty}}}{\alpha^{2}} \ d\alpha \int 
\frac{\Vert f_{x}(x-\alpha)+f_{x}(x+\alpha)-2f_{x}(x)\Vert^{2}_{L^{2}}}{\alpha^{2}} \ d\alpha \right)^{1/2} \\
&\leq& \Vert f \Vert_{\dot H^{1}} \Vert f \Vert_{\dot B^{1/2}_{2,2}} \Vert f \Vert_{\dot H^{3/2}} \\
&\leq& \Vert f \Vert^{2}_{\dot H^{1}} \Vert f \Vert_{\dot H^{3/2}}
\end{eqnarray*}
as well, for $L_{1,2,3}$ we observe that
 \begin{eqnarray*}
L_{1,2,3}&=&-\frac{1}{2}\int \Lambda f \ \int \frac{f(x-\alpha)+f(x+\alpha)-2f(x)}{\alpha^3}    \  \ \int_{0}^{\infty} \delta  e^{-\delta} \int_{0}^{\alpha} f_{x}(x-s)+f_{x}(x+s)-2f_{x}(x)  \ ds \\
&\quad&\times \sin(\frac{\delta}{2} (\Delta_{\alpha}f-\bar\Delta_{\alpha}f))\ d\delta \ d\alpha \ dx \\
 &\leq&\frac{1}{4}\Vert f \Vert_{\dot H^{1}} \int_{0}^{\infty} \delta e^{-\delta} \frac{\Vert f(x-\alpha)+f(x+\alpha)-2f(x) \Vert_{L^{\infty}}}{\vert \alpha \vert^{3}} \\
 &\quad&\times \int_{0}^{\alpha}  {\Vert f_{x}(x-s)-f_{x}(x) \Vert_{L^{2}}}+{\Vert f_{x}(x+s)-f_{x}(x) \Vert_{L^{2}}} \ ds  \ d\alpha \ d\delta \\
 &\leq& \frac{1}{4} \Vert f \Vert_{\dot H^{1}} \int_{0}^{\infty} \delta e^{-\delta}\left( \frac{\Vert f(x-\alpha)-f(x) \Vert_{L^{\infty}}}{\vert \alpha \vert^{3-q-\frac{1}{\bar r}}}  +  \frac{\Vert f_{x}(x+\alpha)-f_{x}(x)\Vert_{L^{\infty}}}{\vert \alpha \vert^{3-q-\frac{1}{\bar r}}}  \right)  \\
 &\quad&\times \left(\int  \frac{\Vert f_{x}(x+s)-f_{x}(x) \Vert^{r} _{L^{2}}}{\vert s \vert^{qr}} \ ds\right)^{1/r}  \ ds\ d\alpha \ d\delta \\
 &\leq& \frac{1}{4} \Vert f \Vert_{\dot H^{1}} \int_{0}^{\infty} \delta e^{-\delta}\left( \frac{\Vert f(x-\alpha)-f(x) \Vert_{L^{\infty}}}{\vert \alpha \vert^{3-q-\frac{1}{\bar r}}}  \right)  \times \left(\int  \frac{\Vert f_{x}(x+s)-f_{x}(x) \Vert^{r} _{L^{2}}}{\vert s \vert^{qr}} \ ds\right)^{1/r}  \ d\alpha \ d\delta \\
  &\leq&  \frac{\Gamma(2)}{4} \Vert f \Vert_{H^{1}} \Vert f \Vert_{\dot B^{2-q-\frac{1}{\bar r}}_{\infty,1}} \Vert f_{x} \Vert_{\dot B^{q-\frac{1}{r}}_{2,r}} \\
\end{eqnarray*}
Then, by choosing  $q=3/4,$ $r=\bar r=2$,  one gets
\begin{eqnarray*}
\vert L_{1,2,3} \vert &\leq& \Vert f \Vert_{H^{1}} \Vert f \Vert_{\dot B^{3/4}_{\infty,1}} \Vert f \Vert_{\dot B^{5/4}_{2,2}}   \\
&\leq& \frac{1}{2} \Vert f \Vert_{\dot H^{1}}  \Vert f \Vert^{1/2}_{\dot B^{1/2}_{\infty,\infty}}\Vert f \Vert^{1/2}_{\dot B^{1}_{\infty,\infty}}  \Vert f \Vert^{1/2}_{\dot B^{3/2}_{2,2}} \Vert f \Vert^{1/2}_{\dot B^{1}_{2,2}} \\
&\leq& \frac{1}{2} \Vert f \Vert^{2}_{\dot H^{1}}  \Vert f \Vert_{\dot H^{3/2}}
\end{eqnarray*}
Since \eqref{10sip} is a dissipative term and by \eqref{a10sip}, we have obtained that
 \begin{equation*} 
  \vert L_{1} \vert \leq \Vert f \Vert^{2}_{H^{1}}   P(\Vert f \Vert_{\dot H^{3/2}}) - \pi \Vert f \Vert^{2}_{\dot H^{1}}+ \pi\frac{K^2}{1+K^2} \Vert f \Vert^{2}_{\dot H^{1}}
 \end{equation*}

 Finally,
  \begin{equation} \label{est1}
 \vert L_{1} \vert \leq \Vert f \Vert^{2}_{H^{1}}   P(\Vert f \Vert_{\dot H^{3/2}})- \frac{\pi}{1+K^2} \Vert f \Vert^{2}_{\dot H^{1}}
 \end{equation}
 where $P(X)=X+X^2.$

   And then integrating in time $s \in [0,T]$ one gets the desired energy inequality.
  Therefore, if $\Vert f_{0} \Vert_{\dot H^{3/2}}$ is smaller than some $C(K)$ that depends only on $K$, then the solution is in $L^{\infty}([0,T],\dot H^{1/2}) \cap L^{2}([0,T], \dot H^{1})$.  This concludes the $\dot H^{1/2}$-estimates.
   \qed
%
%

\begin{thebibliography}{10}

\bibitem{AmbroseST}
D.~{Ambrose}.
\newblock {The zero surface tension limit of two-dimensional interfacial Darcy
  flow.}
\newblock {\em J. Math. Fluid Mech.}, 16:105--143, 2014.

\bibitem{ambrose2004well}
D.M. Ambrose.
\newblock {W}ell-posedness of two-phase {H}ele-{S}haw flow without surface
  tension.
\newblock {\em {E}uropean {J}ournal of {A}pplied {M}athematics},
  15(5):597--607, 2004.

\bibitem{bae2015global}
Hantaek Bae and Rafael Granero-Belinch{\'o}n.
\newblock Global existence for some transport equations with nonlocal velocity.
\newblock {\em Advances in Mathematics}, 269:197--219, 2015.

\bibitem{bahouri2011fourier}
Hajer Bahouri, Jean-Yves Chemin, and Rapha{\"e}l Danchin.
\newblock {\em Fourier analysis and nonlinear partial differential equations},
  volume 343.
\newblock Springer Science \& Business Media, 2011.

\bibitem{bear}
J.~Bear.
\newblock {\em {D}ynamics of fluids in porous media}.
\newblock Dover Publications, 1988.

\bibitem{BCG}
L.C. Berselli, D.~C\'ordoba, and R.~Granero-Belinch\'on.
\newblock {L}ocal solvability and turning for the inhomogeneous {M}uskat
  problem.
\newblock {\em Interfaces and Free Boundaries}, 16(2):175--213, 2014.

\bibitem{berselli2002vanishing}
Luigi~C Berselli.
\newblock Vanishing viscosity limit and long-time behavior for 2d
  quasi-geostrophic equations.
\newblock {\em Indiana University mathematics journal}, 51(4):905--930, 2002.

\bibitem{B}
Oleg~Vladimirovich Besov.
\newblock Investigation of a class of function spaces in connection with
  imbedding and extension theorems.
\newblock {\em Trudy Matematicheskogo Instituta imeni VA Steklova}, 60:42--81,
  1961.

\bibitem{buckley1941mechanism}
SE~Buckley and MC~Leverett.
\newblock Mechanism of fluid displacement in sands.
\newblock {\em Trans. Aime}, 146, 1941.

\bibitem{cameron2017global}
Stephen Cameron.
\newblock Global well-posedness for the 2d Muskat problem with slope less than
  1.
\newblock {\em Analysis \& PDE} 12 (4), 997-1022, 2017.

\bibitem{castro2016mixing}
A~Castro, D~C{\'o}rdoba, and D~Faraco.
\newblock Mixing solutions for the Muskat problem.
\newblock {\em arXiv preprint arXiv:1605.04822}, 2016.

\bibitem{castro2012breakdown}
A.~Castro, D.~Cordoba, C.~Fefferman, and F.~Gancedo.
\newblock {B}reakdown of smoothness for the {M}uskat problem.
\newblock {\em {A}rchive for {R}ational {M}echanics and {A}nalysis},
  208(3):805--909, 2013.

\bibitem{ccfgl}
A.~Castro, D.~Cordoba, C.~Fefferman, F.~Gancedo, and M.~Lopez-Fernandez.
\newblock {R}ayleigh-{T}aylor breakdown for the {M}uskat problem with
  applications to water waves.
\newblock {\em Annals of Math}, 175:909--948, 2012.

\bibitem{ccfgonephase}
Angel Castro, Diego C{\'o}rdoba, Charles Fefferman, and Francisco Gancedo.
\newblock Splash singularities for the one-phase Muskat problem in stable
  regimes.
\newblock {\em Archive for Rational Mechanics and Analysis}, 222(1):213--243,
  2016.

\bibitem{castro2018degraded}
{\'A}ngel Castro, Daniel Faraco, and Francisco Mengual.
\newblock Degraded mixing solutions for the Muskat problem.
\newblock {\em arXiv preprint arXiv:1805.12050}, 2018.

\bibitem{CF}
M.~Cerminara and A.~Fasano.
\newblock {M}odelling the dynamics of a geothermal reservoir fed by gravity
  driven flow through overstanding saturated rocks.
\newblock {\em Journal of Volcanology and Geothermal Research}, 233:37--54,
  2012.

\bibitem{Chae2singularSQG}
Dongho Chae, Peter Constantin, Diego C{\'o}rdoba, Francisco Gancedo, and
  Jiahong Wu.
\newblock Generalized surface quasi-geostrophic equations with singular
  velocities.
\newblock {\em Comm. Pure Appl. Math.}, 65(8):1037--1066, 2012.

\bibitem{chang2016free}
H{\'e}ctor~A Chang-Lara and Nestor Guillen.
\newblock From the free boundary condition for hele-shaw to a fractional
  parabolic equation.
\newblock {\em arXiv preprint arXiv:1605.07591}, 2016.

\bibitem{chen1993hele}
Xinfu Chen.
\newblock The hele-shaw problem and area-preserving curve-shortening motions.
\newblock {\em Archive for rational mechanics and analysis}, 123(2):117--151,
  1993.

\bibitem{CGS}
Ching-Hsiao~Arthur Cheng, Rafael Granero-Belinch\'on, and Steve Shkoller.
\newblock Well-posedness of the Muskat problem with {$H^2$} initial data.
\newblock {\em Advances in Mathematics}, 286:32--104, 2016.

\bibitem{CGSW}
Ching-Hsiao~Arthur Cheng, Rafael Granero-Belinch\'on, Steve Shkoller and Jon Wilkening.
\newblock Rigorous Asymptotic Models of Water Waves.
\newblock {\em Water Waves}, 1--60, 2019.

\bibitem{ccgs-13}
P.~Constantin, D.~Cordoba, F.~Gancedo, Luis Rodr{\'i}guez-Piazza, and R.M.
  Strain.
\newblock On the {M}uskat problem: global in time results in 2d and 3d.
\newblock {\em Amer. J. Math. Vol 138, no.6}, 138(6), 2016.

\bibitem{ccgs-10}
P.~Constantin, D.~Cordoba, F.~Gancedo, and R.M. Strain.
\newblock On the global existence for the {M}uskat problem.
\newblock {\em Journal of the European Mathematical Society}, 15:201--227,
  2013.

\bibitem{constantin1999formation}
P.~Constantin, A.J. Majda, and E.~Tabak.
\newblock {F}ormation of strong fronts in the 2-{D} quasi-geostrophic thermal
  active scalar.
\newblock {\em Nonlinearity}, 7(6):1495, 1994.

\bibitem{constantin1994singular}
P.~Constantin, A.J. Majda, and E.G. Tabak.
\newblock {S}ingular front formation in a model for quasi-geostrophic flow.
\newblock {\em Physics of Fluids}, 6:9, 1994.

\bibitem{Peter}
P.~Constantin and M.~Pugh.
\newblock {G}lobal solutions for small data to the {H}ele-{S}haw problem.
\newblock {\em Nonlinearity}, 6:393--415, 1993.

\bibitem{CGSVfiniteslope}
Peter Constantin, Francisco Gancedo, Roman Shvydkoy, and Vlad Vicol.
\newblock Global regularity for 2d Muskat equations with finite slope.
\newblock {\em To appear in Annales de l'Institut Henri Poincare (C) Non Linear
  Analysis}, 2016.

\bibitem{cor2}
A.~C{\'o}rdoba and D.~C{\'o}rdoba.
\newblock A maximum principle applied to quasi-geostrophic equations.
\newblock {\em Communications in Mathematical Physics}, 249(3):511--528, 2004.

\bibitem{c-c-g10}
A.~Cordoba, D.~C{\'o}rdoba, and F.~Gancedo.
\newblock {I}nterface evolution: the {H}ele-{S}haw and {M}uskat problems.
\newblock {\em Annals of Math}, 173, no. 1:477--542, 2011.

\bibitem{cordoba2009rayleigh}
Antonio Cordoba, Diego Cordoba, and Francisco Gancedo.
\newblock The Rayleigh-Taylor condition for the evolution of irrotational fluid
  interfaces.
\newblock {\em Proceedings of the National Academy of Sciences},
  106(27):10955--10959, 2009.

\bibitem{Cordoba-Cordoba-Gancedo:muskat-3d}
Antonio C{\'o}rdoba, Diego C{\'o}rdoba, and Francisco Gancedo.
\newblock Porous media: the {M}uskat problem in three dimensions.
\newblock {\em Anal. PDE}, 6(2):447--497, 2013.

\bibitem{c-g07}
D.~C{\'o}rdoba and F.~Gancedo.
\newblock {C}ontour dynamics of incompressible 3-{D} fluids in a porous medium
  with different densities.
\newblock {\em Communications in Mathematical Physics}, 273(2):445--471, 2007.

\bibitem{c-g09}
D.~C{\'o}rdoba and F.~Gancedo.
\newblock {A} maximum principle for the {M}uskat problem for fluids with
  different densities.
\newblock {\em Communications in Mathematical Physics}, 286(2):681--696, 2009.

\bibitem{CGO}
D.~C\'ordoba, R.~Granero-Belinch\'on, and R.~Orive.
\newblock {O}n the confined {M}uskat problem: differences with the deep water
  regime.
\newblock {\em Communications in Mathematical Sciences}, 12(3):423--455, 2014.

\bibitem{cponephase}
D.~Cordoba and T.~Pernas-Casta{\~n}o.
\newblock {N}on-splat singularity for the one-phase {M}uskat problem.
\newblock {\em Transaction of the American Mathematical Society},
  369(1):711--754, 2017.

\bibitem{cordoba2010absence}
Diego C{\'o}rdoba and Francisco Gancedo.
\newblock Absence of squirt singularities for the multi-phase Muskat problem.
\newblock {\em Communications in Mathematical Physics}, 299(2):561--575, 2010.

\bibitem{cordoba2015note}
Diego C{\'o}rdoba, Javier G{\'o}mez-Serrano, and Andrej Zlato{\v{s}}.
\newblock A note on stability shifting for the Muskat problem.
\newblock {\em Phil. Trans. R. Soc. A}, 373(2050):20140278, 2015.

\bibitem{cordoba2017note}
Diego C{\'o}rdoba, Javier G{\'o}mez-Serrano, and Andrej Zlato{\v{s}}.
\newblock A note on stability shifting for the Muskat problem, ii: From stable
  to unstable and back to stable.
\newblock {\em Analysis \& PDE}, 10(2):367--378, 2017.

\bibitem{CL}
Diego Cordoba and Omar Lazar.
\newblock Global well-posedness for the 2d stable Muskat problem in $H^{3/2}$.
\newblock {\em arXiv preprint arXiv:1803.07528}, 2018.

\bibitem{coutand2007well}
Daniel Coutand and Steve Shkoller.
\newblock Well-posedness of the free-surface incompressible Euler equations
  with or without surface tension.
\newblock {\em Journal of the American Mathematical Society}, 20(3):829--930,
  2007.

\bibitem{coutand2016impossibility}
Daniel Coutand and Steve Shkoller.
\newblock On the impossibility of finite-time splash singularities for vortex
  sheets.
\newblock {\em Archive for Rational Mechanics and Analysis}, 221(2):987--1033,
  2016.

\bibitem{Darcy}
Henry Darcy.
\newblock {\em Les fontaines publiques de la ville de Dijon: exposition et
  application}
\newblock Victor Dalmont, 1856.

\bibitem{DMP}
L~Dawson, H~McGahagan, and G~Ponce.
\newblock On the decay properties of solutions to a class of Schr{\"o}dinger
  equations.
\newblock {\em Proceedings of the American Mathematical Society},
  136(6):2081--2090, 2008.

\bibitem{elliott1982weak}
Charles~M Elliott and John~R Ockendon.
\newblock {\em Weak and variational methods for moving boundary problems},
  volume~59.
\newblock Pitman Publishing, 1982.

\bibitem{escher2011generalized}
Joachim Escher, Anca-Voichita Matioc, and Bogdan-Vasile Matioc.
\newblock A generalized {R}ayleigh-{T}aylor condition for the {M}uskat problem.
\newblock {\em Nonlinearity}, 25(1):73--92, 2012.

\bibitem{e-m10}
Joachim Escher and Bogdan-Vasile Matioc.
\newblock On the parabolicity of the {M}uskat problem: Well-posedness,
  fingering, and stability results.
\newblock {\em Zeitschrift f{\"u}r Analysis und ihre Anwendungen},
  30(2):193--218, 2011.

\bibitem{emw15}
Joachim Escher, Bogdan-Vasile Matioc, and Christoph Walker.
\newblock The domain of parabolicity for the {M}uskat problem.
\newblock {\em Indiana Univ. Math. J.}, 67(2):679--737, 2018.

\bibitem{escher1998center}
Joachim Escher and Gieri Simonett.
\newblock A center manifold analysis for the mullins--sekerka model.
\newblock {\em journal of differential equations}, 143(2):267--292, 1998.

\bibitem{escher1997classical}
Joachim Escher, Gieri Simonett, et~al.
\newblock Classical solutions for hele-shaw models with surface tension.
\newblock {\em Advances in Differential Equations}, 2(4):619--642, 1997.

\bibitem{FIL}
Charles Fefferman, Alexandru~D Ionescu, Victor Lie, et~al.
\newblock On the absence of splash singularities in the case of two-fluid
  interfaces.
\newblock {\em Duke Mathematical Journal}, 165(3):417--462, 2016.

\bibitem{forster2017piecewise}
Clemens F{\"o}rster and L{\'a}szl{\'o} Sz{\'e}kelyhidi~Jr.
\newblock Piecewise constant subsolutions for the Muskat problem.
\newblock {\em Communications in Mathematical Physics} 363(3),1051--1080, 2018.

\bibitem{Friedlander3}
Susan Friedlander, Walter Rusin, and Vlad Vicol.
\newblock On the supercritically diffusive magnetogeostrophic equations.
\newblock {\em Nonlinearity}, 25(11):3071, 2012.

\bibitem{Friedlander5}
Susan Friedlander, Walter Rusin, Vlad Vicol, and AI~Nazarov.
\newblock The magneto-geostrophic equations: a survey.
\newblock {\em Proc. of the St. Petersburg Mathematical Society, Volume XV:
  Advances in Mathematical Analysis of Partial Differential Equations}, 2014.

\bibitem{Friedlander2}
Susan Friedlander and Vlad Vicol.
\newblock Global well-posedness for an advection--diffusion equation arising in
  magneto-geostrophic dynamics.
\newblock In {\em Annales de l'Institut Henri Poincare (C) Non Linear
  Analysis}, volume~28, pages 283--301. Elsevier, 2011.

\bibitem{Friedlander4}
Susan Friedlander and Vlad Vicol.
\newblock On the ill/well-posedness and nonlinear instability of the
  magneto-geostrophic equations.
\newblock {\em Nonlinearity}, 24(11):3019, 2011.

\bibitem{F}
A.~Friedman.
\newblock {F}ree boundary problems arising in tumor models.
\newblock {\em Atti Accad. Naz. Lincei Cl. Sci. Fis. Mat. Natur. Rend.
  Lincei,}, 9(3-4), 2004.

\bibitem{gancedo2008existence}
F.~Gancedo.
\newblock {E}xistence for the {$\alpha$}-patch model and the {QG} sharp front
  in {S}obolev spaces.
\newblock {\em Advances in Mathematics}, 217(6):2569--2598, 2008.

\bibitem{gancedo2017muskat}
Francisco Gancedo, Eduardo Garcia-Juarez, Neel Patel, and Robert~M Strain.
\newblock On the muskat problem with viscosity jump: Global in time results.
\newblock {\em Advances in Mathematics}, 345, 552-597, 2019.

\bibitem{gancedo2014absence}
Francisco Gancedo and Robert~M Strain.
\newblock Absence of splash singularities for surface quasi-geostrophic sharp
  fronts and the {M}uskat problem.
\newblock {\em Proceedings of the National Academy of Sciences},
  111(2):635--639, 2014.

\bibitem{GG}
J.~G\'omez-Serrano and R.~Granero-Belinch\'on.
\newblock On turning waves for the inhomogeneous {M}uskat problem: a
  computer-assisted proof.
\newblock {\em Nonlinearity}, 27(6):1471--1498., 2014.

\bibitem{G}
R.~Granero-Belinch{\'o}n.
\newblock Global existence for the confined {M}uskat problem.
\newblock {\em SIAM Journal on Mathematical Analysis}, 46(2):1651--1680, 2014.

\bibitem{GScr}
R.~Granero-Belinch{\'o}n \& S. Scrobogna,
\newblock Asymptotic models for free boundary flow in porous media,
\newblock {\em Physica D: Nonlinear Phenomena}, 2019.

\bibitem{GS}
R~Granero-Belinch{\'o}n and S~Shkoller.
\newblock Well-posedness and decay to equilibrium for the muskat problem with
  discontinuous permeability,(2016). preprint.
\newblock {\em To appear in Transactions of the American Mathematical Society,
  arXiv preprint arXiv:1611.06147}.

\bibitem{granero2013inhomogeneous}
Rafael Granero~Belinch{\'o}n et~al.
\newblock The inhomogeneous Muskat problem.
\newblock 2013.

\bibitem{hassanizadeh1990mechanics}
S~Majid Hassanizadeh and William~G Gray.
\newblock Mechanics and thermodynamics of multiphase flow in porous media
  including interphase boundaries.
\newblock {\em Advances in water resources}, 13(4):169--186, 1990.

\bibitem{H-S}
H.~S. Hele-Shaw.
\newblock The flow of water.
\newblock {\em Nature}, 58:34--36, 1898.

\bibitem{HeleShaw:motion-viscous-fluid-parallel-plates}
H.~S. Hele-Shaw.
\newblock On the motion of a viscous fluid between two parallel plates.
\newblock {\em Trans. Royal Inst. Nav. Archit.}, 40:218, 1898.

\bibitem{hornung1997homogenization}
U.~Hornung.
\newblock {\em {H}omogenization and porous media}, volume~6.
\newblock Springer Verlag, 1997.

\bibitem{KiselevSQG}
A.~Kiselev, F.~Nazarov, and A.~Volberg.
\newblock Global well-posedness for the critical 2{D} dissipative
  quasi-geostrophic equation.
\newblock {\em Invent. Math.}, 167(3):445--453, 2007.

\bibitem{Omar}
Omar Lazar.
\newblock Global existence for the critical dissipative surface
  quasi-geostrophic equation.
\newblock {\em Comm. Math. Phys.}, 322(1):73--93, 2013.

\bibitem{lemarie2016navier}
Pierre~Gilles Lemari{\'e}-Rieusset.
\newblock {\em The Navier-Stokes problem in the 21st century}.
\newblock Chapman and Hall/CRC, 2016.

\bibitem{majda1996two}
A.J. Majda and E.G. Tabak.
\newblock {A} two-dimensional model for quasigeostrophic flow: comparison with
  the two-dimensional {E}uler flow.
\newblock {\em Physica D: Nonlinear Phenomena}, 98(2-4):515--522, 1996.

\bibitem{majda2002vorticity}
Andrew~J Majda and Andrea~L Bertozzi.
\newblock {\em Vorticity and incompressible flow}, volume~27.
\newblock Cambridge University Press, 2002.

\bibitem{matioc2017well}
Anca-Voichita Matioc and Bogdan-Vasile Matioc.
\newblock Well-posedness and stability results for a quasilinear periodic
  muskat problem.
\newblock {\em Journal of Differential Equations}, 266(9), 5500-5531, 2019.

\bibitem{matioc2016muskat}
Bogdan-Vasile Matioc.
\newblock The muskat problem in 2d: equivalence of formulations,
  well-posedness, and regularity results.
\newblock {\em Analysis \& PDE}, 12(2), 281-332, 2018.

\bibitem{MR3841857}
Bogdan-Vasile Matioc.
\newblock Viscous displacement in porous media: the Muskat problem in 2{D}.
\newblock {\em Trans. Amer. Math. Soc.}, 370(10):7511--7556, 2018.

\bibitem{matioc2018well}
Bogdan-Vasile Matioc.
\newblock Well-posedness and stability results for some periodic Muskat
  problems.
\newblock {\em arXiv preprint arXiv:1804.10403}, 2018.

\bibitem{Moffatt}
HK~Moffatt and DE~Loper.
\newblock The magnetostrophic rise of a buoyant parcel in the earth's core.
\newblock {\em Geophysical Journal International}, 117(2):394--402, 1994.

\bibitem{Muskat}
M.~Muskat.
\newblock {T}he flow of homogeneous fluids through porous media.
\newblock {\em Soil Science}, 46(2):169, 1938.

\bibitem{Musk}
Morris Muskat.
\newblock Two fluid systems in porous media. the encroachment of water into an
  oil sand.
\newblock {\em Physics}, 5(9):250--264, 1934.

\bibitem{Muskat:porous-media}
Morris Muskat.
\newblock The flow of fluids through porous media.
\newblock {\em Journal of Applied Physics}, 8(4):274--282, 1937.

\bibitem{NB}
D.A. Nield and A.~Bejan.
\newblock {\em {C}onvection in porous media}.
\newblock Springer Verlag, 2006.

\bibitem{otto1999evolution}
Felix Otto.
\newblock Evolution of microstructure in unstable porous media flow: a
  relaxational approach.
\newblock {\em Communications on Pure and Applied Mathematics: A Journal Issued
  by the Courant Institute of Mathematical Sciences}, 52(7):873--915, 1999.

\bibitem{otto2001evolution}
Felix Otto.
\newblock Evolution of microstructure: an example.
\newblock In {\em Ergodic theory, analysis, and efficient simulation of
  dynamical systems}, pages 501--522. Springer, 2001.

\bibitem{patel2017large}
Neel Patel and Robert~M Strain.
\newblock Large time decay estimates for the Muskat equation.
\newblock {\em Communications in Partial Differential Equations},
  42(6):977--999, 2017.

\bibitem{pernas2017local}
Tania Pernas-Casta{\~n}o.
\newblock Local-existence for the inhomogeneous Muskat problem.
\newblock {\em Nonlinearity}, 30(5):2063, 2017.

\bibitem{pruess2016muskat}
Jan Pruess and Gieri Simonett.
\newblock On the Muskat flow.
\newblock {\em Evolution Equations and Control Theory}, 5:631--645, 2016.

\bibitem{Rayleigh:instability-jets}
Lord Rayleigh.
\newblock On the instability of jets.
\newblock {\em Proceedings of the London Mathematical Society}, s1-10(1):4--13,
  1878.

\bibitem{RodrigoQG}
Jos{\'e}~Luis Rodrigo.
\newblock On the evolution of sharp fronts for the quasi-geostrophic equation.
\newblock {\em Comm. Pure Appl. Math.}, 58(6):821--866, 2005.

\bibitem{RS}
Thomas Runst and Winfried Sickel.
\newblock {\em Sobolev spaces of fractional order, Nemytskij operators, and
  nonlinear partial differential equations}, volume~3.
\newblock Walter de Gruyter, 1996.

\bibitem{SaffTay}
P.~G. Saffman and Geoffrey Taylor.
\newblock The penetration of a fluid into a porous medium or {H}ele-{S}haw cell
  containing a more viscous liquid.
\newblock {\em Proc. Roy. Soc. London. Ser. A}, 245:312--329. (2 plates), 1958.

\bibitem{SCH}
M.~Siegel, R.E. Caflisch, and S.~Howison.
\newblock {G}lobal existence, singular solutions, and ill-posedness for the
  {M}uskat problem.
\newblock {\em Communications on {P}ure and {A}pplied {M}athematics},
  57(10):1374--1411, 2004.

\bibitem{szekelyhidi2012relaxation}
L{\'a}szl{\'o} Sz{\'e}kelyhidi~Jr.
\newblock Relaxation of the incompressible porous media equation.
\newblock {\em Ann. Sci. {\'E}c. Norm. Sup{\'e}r.(4)}, 45(3):491--509, 2012.

\bibitem{Tartar:incompressible-porous-medium-homogenization}
Luc Tartar.
\newblock Incompressible fluid flow in a porous medium-convergence of the
  homogenization process.
\newblock In {\em Nonhomogeneous media and vibration theory}. 1980.
\newblock E. S\'anchez-Palencia.

\bibitem{thornton2014hele}
Anthony~R Thornton, Avraham~J van~der Horn, Elena Gagarina, Wout Zweers,
  Devaraj van~der Meer, and Onno Bokhove.
\newblock Hele-shaw beach creation by breaking waves: a mathematics-inspired
  experiment.
\newblock {\em Environmental fluid mechanics}, 14(5):1123--1145, 2014.

\bibitem{tofts2017existence}
Spencer Tofts.
\newblock On the existence of solutions to the Muskat problem with surface
  tension.
\newblock {\em Journal of Mathematical Fluid Mechanics}, 19(4):581--611, 2017.

\end{thebibliography}

\end{document}